\documentclass[letterpaper]{article} % DO NOT CHANGE THIS
\usepackage{aaai24}  % DO NOT CHANGE THIS
\usepackage{times}  % DO NOT CHANGE THIS
\usepackage{helvet}  % DO NOT CHANGE THIS
\usepackage{courier}  % DO NOT CHANGE THIS
\usepackage[hyphens]{url}  % DO NOT CHANGE THIS
\usepackage{graphicx} % DO NOT CHANGE THIS
\usepackage{natbib}  % DO NOT CHANGE THIS AND DO NOT ADD ANY OPTIONS TO IT
\usepackage{caption} % DO NOT CHANGE THIS AND DO NOT ADD ANY OPTIONS TO IT
\usepackage{algorithm}
\usepackage{algorithmic}

\usepackage{newfloat}
\usepackage{listings}
\DeclareCaptionStyle{ruled}{labelfont=normalfont,labelsep=colon,strut=off} % DO NOT CHANGE THIS
\floatstyle{ruled}
\newfloat{listing}{tb}{lst}{}
\floatname{listing}{Listing}
\usepackage{amsmath}
\usepackage{mathtools}
\usepackage{mathrsfs}
\usepackage{amssymb}
\usepackage{booktabs}
\usepackage{tablefootnote}
\usepackage{pifont}% http://ctan.org/pkg/pifont
\usepackage[shortlabels]{enumitem}
\setlist[enumerate]{noitemsep, topsep=0pt}
\usepackage{adjustbox}
\usepackage{xspace}
\usepackage{multirow}
\usepackage{amsthm}

\newtheorem{theorem}{Theorem}
\newtheorem{lemma}{Lemma}
\newtheorem{proposition}{Proposition}
\newtheorem{definition}{Definition}
\newtheorem{assumption}{Assumption}
\newtheorem{remark}{Remark}
\newtheorem{corollary}{Corollary}
\usepackage{color}

\title{A Bregman Proximal Stochastic Gradient Method with Extrapolation for Nonconvex Nonsmooth Problems}
\author{
	Qingsong Wang\textsuperscript{\rm 1},
	Zehui Liu\textsuperscript{\rm 2}, 
	Chunfeng Cui\textsuperscript{\rm 2}\thanks{Corresponding author.}, 
	Deren Han\textsuperscript{\rm 2}
}
\affiliations{
	\textsuperscript{\rm 1}School of Mathematics and Computational Science, Xiangtan University, China\\
	\textsuperscript{\rm 2}LMIB, School of Mathematical Sciences, Beihang University, China\\
	nothing2wang@hotmail.com, \{liuzehui, chunfengcui, handr\}@buaa.edu.cn
}

\begin{document}
	
	\maketitle
	
	\begin{abstract}
		In this paper, we explore a specific optimization problem that involves the combination of a differentiable nonconvex function and a nondifferentiable function. The differentiable component lacks a global Lipschitz continuous gradient, posing challenges for optimization.  To address this issue and accelerate the convergence, we propose a Bregman proximal stochastic gradient method with extrapolation (BPSGE), which only requires smooth adaptivity of the differentiable part.  Under the variance reduction framework, we not only analyze the subsequential and global convergence of the proposed algorithm under certain conditions,  but also analyze the sublinear convergence rate of the subsequence, and the complexity of the algorithm, revealing that the BPSGE algorithm requires at most  $\mathcal{O}(\varepsilon^{-2})$ iterations in expectation to attain an $\varepsilon$-stationary point.  To validate the effectiveness of our proposed algorithm, we conduct numerical experiments on three real-world applications: graph regularized nonnegative matrix factorization (NMF), matrix factorization with weakly-convex regularization, and NMF with nonconvex sparsity constraints. These experiments demonstrate that BPSGE is faster than the baselines without extrapolation.
	\end{abstract}
	
	\section{Introduction}
	In this paper, we  consider the following nonsmooth nonconvex optimization problem
	\begin{eqnarray}
		\min_{x\in \overline{C}}\quad \Phi(x):=f(x)+h(x), \label{target_function}
	\end{eqnarray}
	where $\overline{C}$ denotes the closure of $C$, which is a nonempty, convex, and open set in $\mathbb{R}^{d}$, $f$ is a continuously differentiable  (may be nonconvex) function which can be written as $f(x)=\frac{1}{n}\sum_{i=1}^{n}f_{i}(x)$, and $h$ is an extended valued function (maybe nonconvex), which is a regularizer promoting low-complexity structures such as sparsity \cite{Donoho95, FanL01, Zhang10} or non-negativity \cite{LeeS99, HeXZZC11} in the solution. Throughout this paper, the usual restrictive requirement of global Lipschitz continuity of the gradient of $f$ is not needed  \cite{Lan2020First, Gillis20}. Many applications can be categorized in the optimization problem \eqref{target_function}, such as matrix and tensor factorization \cite{ComonGLM08, KoldaB09, CheW20}, supervised neural network model \cite{HasannasabHNPSS2020},  and Poisson deconvolution/inverse problems \cite{BolteSTV18First}. 
	
	\textbf{Deterministic Bregman type methods.} The Bregman proximal gradient (BPG) algorithm \cite{BolteSTV18First} is a state-of-the-art method for addressing optimization problems characterized by the absence of global  Lipschitz continuous gradients. It is simple and far-reaching, 
	%and based on adapting the geometry to the objective through the Bregman distance paradigm. Specifically, this paradigm involves substituting 
	which substitutes the customary upper quadratic approximation of a smooth function with a more comprehensive proximity measure \cite{BauschkeBT17}.  Zhang et al. \cite{ZhangBM0C19} proposed an extrapolation version of the BPG algorithm, denoted by BPGE, under the assumption that  $h(x)$ is convex. Recently, Mukkamala et al. \cite{MukkamalaOPS20} introduced an inertial variant of the BPG method, referred to as the CoCaIn BPG algorithm, which tunes the inertial parameter by convex-concave backtracking.
	%This backtracking strategy capitalizes on the ability to construct both lower and upper approximations for smooth adaptable functions. 
	Compared with the original BPG method \cite{BolteSTV18First}, the CoCaIn BPG algorithm exhibits superior numerical performance. 
	Moreover, the BPG framework facilitates the discovery of novel techniques tailored to specific problem domains. For instance,  Teboulle and Vaisbourd \cite{TeboulleV20} explored new decomposition settings of the nonnegative matrix factorization (NMF) problem with sparsity constraints. Additionally, an alternative perspective is to consider the BPG-type algorithm by the splitting method. Ahookhos et al. \cite{AhookhoshTP21} proposed a Bregman forward-backward splitting line search method, which demonstrates locally superlinear convergence to nonisolated local minima by leveraging the Bregman forward-backward envelope \cite{BauschkeDL18, LaudeOC20}. Furthermore, Wang et al. \cite{WangTOW22} investigated a Bregman and inertial extension of the forward-reflected-backward algorithm \cite{MalitskyT20} under relative smoothness conditions. %This extension involves incorporating inertial steps in the dual space and allows for the consideration of potentially negative inertial values. 
	Numerous other works have explored the Bregman gradient method framework to tackle the absence of global  Lipschitz continuous gradients, including  \cite{ReemRP19, ZhaoDRW22, ZhuDLZ21, DragomirTdB22}.

	\begin{table*}[!ht]	
		\centering
		\begin{adjustbox}{max width=\textwidth}
			\begin{tabular}{c|c| c| c |c |c}\hline
				Versions& Algorithm&$h(x)$&Inertial&Sequence convergence &  Complexity\cr\hline
				\multirow{6}{*}{Deterministic} & BPG \cite{BolteSTV18First}& nonconvex & no &subsequential/global &-\\\cline{2-6}
				&BPGE \cite{ZhangBM0C19}& convex &yes & subsequential/global &- \\\cline{2-6}
				& CoCaIn BPG \cite{MukkamalaOPS20}& weakly-convex &yes & subsequential/global &- \\\cline{2-6}
				& BBPG \cite{TeboulleV20} & nonconvex & no &global& -\\\cline{2-6}
				& $i^{*}$FRB \cite{WangTOW22} & nonconvex &no & global&- \\\cline{2-6}
				& BELLA \cite{AhookhoshTP21}& nonconvex &no & subsequential/global&-\\\hline
				\multirow{5}{*}{Stochastic} & SMD \cite{ZhangH18} & convex & no &- & $\mathcal{O}(\varepsilon^{-2})$\\\cline{2-6}
				& SVRAMD \cite{LiWZC22} & convex &no &  - & $\mathcal{O}(n\varepsilon^{-\frac{2}{3}}+\varepsilon^{-\frac{5}{3}})$\\\cline{2-6}
				&BFinito \cite{LatafatTAP22} & nonconvex  & no & subsequential/global &-\\\cline{2-6}
				& BPSG \cite{WangH23} & nonconvex & no &subsequential/global  & - \\\cline{2-6}
				& BPSGE (Algorithm \ref{BPSGE}) & weakly-convex & yes &subsequential/global  & $\mathcal{O}(\varepsilon^{-2})$\\\hline
			\end{tabular}
		\end{adjustbox}
		\caption{Summary of the properties of BPSGE (Algorithm \ref{BPSGE}) and several state-of-the-art methods. ``Complexity'' means the complexity (in expectation) to obtain an $\varepsilon$-stationary point (Definition \ref{stationary-point}) of $\Phi$ and ``-'' means not given.}
		\label{results}
	\end{table*}

	\textbf{Stochastic Bregman type methods.} For large-scale datasets, the cost associated with computing the full gradient can become prohibitively high.  To address this concern, the full gradient can be replaced with a stochastic gradient estimator \cite{RobbinsM1951, Bottou10} which has yielded remarkable achievements in the field of machine learning \cite{LinLF2020book, Lan2020First}. For the optimization problem without the global  Lipschitz continuous gradients, Zhang and He \cite{ZhangH18} conducted a study on the non-asymptotic stationary convergence behaviour of stochastic mirror descent under certain conditions.  In a similar vein, Li et al. \cite{LiWZC22} proposed a simple and generalized algorithmic framework for incorporating variance reduction into adaptive mirror descent algorithms. However, both approaches \cite{ZhangH18, LiWZC22} require $h(\cdot)$ to be a convex function. To address the limitations posed by convexity assumptions, Latafat et al. \cite{LatafatTAP22} introduced a  Bregman incremental aggregated method that extends the Finito/MISO techniques\cite{DefazioDC14, Mairal15} to non-Lipschitz and nonconvex scenarios. However, it is noteworthy that this approach is memory-intensive and demands periodic computation of the full gradient, which is expensive for large-scale problems. Furthermore, the analysis carried out by Latafat et al. \cite{LatafatTAP22} only encompasses two specific variance reduction stochastic estimators, raising doubts about the generalizability of the convergence results to other estimators. To address the aforementioned concerns,  Wang and Han \cite{WangH23} introduced a stochastic version of the BPG algorithm \cite{BolteSTV18First}, known as the Bregman proximal stochastic gradient (BPSG) method which does not assume $h(\cdot)$ is convex. Under a general framework of variance reduction, they also established the convergence properties of the generated sequence in terms of subsequential and global convergence. 
	
	A summary of the aforementioned algorithms is presented in Table \ref{results}. 
	
	Notwithstanding the notable advancements, there remain areas that warrant further improvement.  Firstly, 
	%the complexity  analysis for the proposed algorithms \cite{LatafatTAP22, WangH23} has not been addressed when $h(\cdot)$ is nonconvex.  
	when $h(\cdot)$ is nonconvex, existing stochastic methods such as BFinito \cite{LatafatTAP22} and BPSG \cite{WangH23}   only analyzed the subsequence and global convergence of the proposed algorithms, yet the sublinear convergence rate of the subsequence and the complexity of the algorithms are unknown. 
	Secondly, the accelerated versions of the Bregman stochastic gradient methods  %in \cite{ZhangH18, LiWZC22, LatafatTAP22, WangH23} 
	have not been taken into account. 
	%From numerical experiments   
	Numerically, the incorporation of accelerated techniques, such as the heavy ball \cite{Polyak64} and the Nesterov acceleration \cite{Nesterov1983}, with BPSG can further accelerate the numerical performance \cite{LinLF2020book}.

	In this paper, we introduce the Bregman proximal stochastic gradient method with extrapolation (BPSGE) to solve the nonconvex nonsmooth optimization problem \eqref{target_function}. Our main contributions addressed in this article are as follows:
	\begin{itemize}
		\item With a general variance reduction framework (see Definition \ref{vr_definition}), we establish the sublinear convergence rate for the subsequence %subsequential sequence 
		generated by the proposed algorithm.
		\item Under the assumption of $C=\mathbb{R}^{d}$, we show that the BPSGE algorithm  requires at most $\mathcal{O}(\varepsilon^{-2})$ iterations in expectation to attain an $\varepsilon$-stationary point. Additionally, we establish the global convergence of the {sequence} generated by   BPSGE, leveraging the Kurdyka-{\L}ojasiewicz (K{\L}) property. 
		\item Numerical experiments are conducted on three distinct problems: graph regularized NMF, matrix factorization (MF) with weakly-convex regularization, and  NMF with nonconvex sparsity constraints. The results of these experiments highlight the favourable performance and enhanced efficiency of the BPSGE algorithm when compared to corresponding algorithms without extrapolation.
	\end{itemize}
	
	The rest of this paper is organized as follows. Section \ref{preliminary} provides some relevant deﬁnitions and results. We present a detailed formulation of the BPSGE algorithm and prove its convergence and convergence rate results in Section \ref{algorithm} and Section \ref{convergence_analysis}, respectively. In Section \ref{numercial_experiments} we use three applications to compare BPSGE with several other algorithms. Finally, we draw conclusions in Section \ref{conclusion}.

	\section{Preliminary}\label{preliminary}
	%In this section, we summarize some useful definitions.  
	%%%%%%%%%%%%%%%%%%
	
	%%%%%%%%%%%%%%%%%%%%%%%%%%%%%%%%%%
	\begin{definition}\label{def:kernel}
		(\cite{AuslenderT06,BolteSTV18First} kernel generating distance). Let $C$ be a nonempty, convex, and open subset of $\mathbb{R}^{d}$. Associated with $C$, a function $\psi:\mathbb{R}^{d} \rightarrow (-\infty, +\infty]$ is called a kernel generating distance if it satisfies the following conditions:
		\begin{itemize}
			\item $\psi$ is proper, lower semicontinuous, and convex, with $\text{dom } \psi$ $\subset$ $\overline{C}$ and  $\text{dom } \partial \psi$  $=C$.
			\item  $\psi$ is $C^{1}$ on $\text{int dom } \psi \equiv C$.
		\end{itemize}
		The class of kernel-generating distances is denoted by $\mathcal{G}(C)$. Given $\psi\in\mathcal{G}(C)$, we  define the proximity measure $D_{\psi} : \text{dom } \psi \times \text{int dom }\psi \rightarrow \mathbb{R}_{+}$ by
		\begin{eqnarray}
			D_{\psi}(x, y) := \psi (x) - \psi (y) - \langle \nabla\psi(y), x - y\rangle.
		\end{eqnarray}
		The proximity measures $D_{\psi}$ is called Bregman distance \cite{Bregman67The}. It measures the proximity of $x$ and $y$. 
	\end{definition}
	Indeed, $\psi$ is convex if and only if $D_{\psi}(x,y) \ge 0$ for any $x\in\text{dom }\psi$ and $y\in\text{int dom }\psi$.
	%due to the gradient inequality.}
	
	\begin{definition}\label{Ll-smooth}
		(\cite{MukkamalaOPS20} $(\bar{L},\underline{L})$-smooth adaptable) Given $\psi\in\mathcal{G}(C)$, let $f:\mathcal{X}\rightarrow(-\infty,+\infty]$ be a proper and lower semi-continuous function with $\mathrm{dom}\,\psi\subset\mathrm{dom}\,f$, which is continuously differentiable on $C$. We say $(f, \psi)$ is $(\bar{L},\underline{L})$- smooth adaptable  on $C$ if there exist $\bar{L}>0$ and $\underline{L}\ge0$ such that for any $x,y\in C$,
		\begin{eqnarray}
			f(x)-f(y)-\langle\nabla f(y),x-y\rangle\le \bar{L} D_{\psi}(x,y),\label{L_upper}
		\end{eqnarray}
		and
		\begin{eqnarray}
			-\underline{L}D_{\psi}(x,y)\le f(x)-f(y)-\langle \nabla f(y),x-y\rangle.\label{L_lower}
		\end{eqnarray}
	\end{definition}
	%Whenever $f$ is convex, 
	If $\underline{L}=\bar{L}$, it recovers \cite[Definition 2.2]{BolteSTV18First}. Suppose  $f$ is convex. If $\underline{L}=0$, this definition recovers \cite[Lemma 1]{BauschkeBT17} and \cite[Definition 1.1]{LuFN18}. 
	
	\begin{definition} \label{stationary-point} 
		(\cite{Lan2020First} $\epsilon$-stationary point) 
		Given $\epsilon>0$, a point $x^{*}$ is called an $\epsilon$-stationary point of function $\Phi(x)$ if
		\[
		\mbox{dist}(0,\partial \Phi(x^{*}))\le\epsilon.
		\]
	\end{definition}
	
	%%%%%%%%%%%%%%%
	\section{Algorithm}\label{algorithm} 
	The classic BPG algorithm \cite{BolteSTV18First}  is given by
	\begin{align*}
		x_{k+1}\in\underset{x\in\bar{C}}{\text{argmin}} \,\, h(x)+\langle \nabla f(x_{k}),x-x_{k}\rangle +\frac{1}{\eta_{k}}D_{\psi}(x,x_{k}),
	\end{align*}
	where $\eta_{k}>0$ is the stepsize. If $\psi=\frac{1}{2}\|\cdot\|^{2}$, it reduces to  the   proximal gradient method. 
	When $n$ is large, replacing the full gradient $\nabla f(x_{k})$   by the stochastic gradient $\tilde{\nabla}f(x_{k})$ can significantly reduce the computational cost, as shown in the BPSG algorithm \cite{WangH23}.
	
	In this section, we introduce the BPSGE algorithm, i.e., an extrapolation or acceleration version of the BPSG algorithm \cite{WangH23},   to solve the nonconvex nonsmooth optimization problem \eqref{target_function}.  
	Before presenting the algorithm framework of BPSGE, we make the following assumptions.
	\begin{assumption} \label{assume_01}
		We assume   the following three conditions hold:
		\begin{itemize}
			\item $\psi\in\mathcal{G}(C)$ is a kernel generating distance given by Definition~\ref{def:kernel} with $\overline{C}=\overline{\text{dom }h}$.
			\item $h:\mathbb{R}^{d}\rightarrow(-\infty,+\infty]$ is a proper and lower semicontinuous function with $\text{dom }h\cap C\neq\emptyset$.
			\item $f:\mathbb{R}^{d}\rightarrow(-\infty,+\infty]$ is a proper and lower semicontinuous function with $\text{dom }\psi\subset\text{dom }f$, and  $\psi$ is continuously differentiable on $C$.
			%\item[(iv)] $\mathcal{V}(\Phi):=\inf\{\Phi(x):\,\,x\in\overline{C}\}>-\infty$.
		\end{itemize}
	\end{assumption}
	%Based  on the assumption described above, 
	Assumption~\ref{assume_01} is quite weak. In addition, we make the following assumptions.
	%%%%%%%%%%
	\begin{assumption}\label{assume_02}
		\begin{itemize}
			\item The function $\psi:\mathbb{R}^{d}\rightarrow(-\infty,+\infty]$ is $\sigma$-strongly convex on $C$. Let $\sigma=1$ for simplicity.
			\item There exists $\alpha\in\mathbb{R}_{+}$ such that $h(\cdot)+\frac{\alpha}{2}\|\cdot\|^{2}$ is convex.
			\item The pair of functions $(f,\psi)$ is $(\bar{L},\underline{L})$-smooth adaptable on $C$. 
		\end{itemize}
	\end{assumption} 
	
	%Based on these two assumptions, 
	Now we describe the framework of the BPSGE algorithm as follows.
	
	\begin{algorithm}[t]
		\caption{BPSGE: Bregman proximal stochastic gradient method with extrapolation}
		\label{BPSGE}
		{\bfseries Input:}  Choose $\psi\in\mathcal{G}(C)$ with $C\equiv\text{int dom }\psi$ such that $(f,\psi)$ is $(\bar{L},\underline{L})$-smooth adaptable  on $C$.\\
		{\bfseries Initialization:} $x_{-1}=x_{0}\in\mathrm{int}\,\mathrm{dom}\,\psi$, $k_{\max}$, two constants $\delta,\epsilon$ such that $0<\epsilon<\delta<1$. \\
		{\bfseries Update:} 
		For $k=0,1,\dots,k_{\max}$,   
		\begin{algorithmic}[1]
			\STATE Compute an extrapolation parameter $\beta_{k}\in[0,1)$ such that
			\begin{eqnarray}
				D_{\psi}(x_{k},\bar{x}_{k})\le \frac{\delta-\epsilon} {1+\underline{L}\eta_{k-1}}D_{\psi}(x_{k-1},x_{k}), \label{extra_ineq}
			\end{eqnarray}
			where  $\bar{x}_{k}=x_{k}+\beta_{k}(x_{k}-x_{k-1})\in\mathrm{int}\,\mathrm{dom}\, \psi$.
			\STATE Compute the stochastic gradient $\tilde{\nabla}f(\bar{x}_{k})$ with the minibatch  $B_{k}$. 
			\STATE Set $0<\eta_{k}\le\min\{\eta_{k-1},\bar{L}^{-1}\}$ and compute $x_{k+1}$ by 
			\begin{eqnarray}
				\begin{aligned}
					x_{k+1}\in\underset{x\in\bar{C}}{\text{argmin}}\,\, &h(x)+\langle \tilde{\nabla}f(\bar{x}_{k}),x-\bar{x}_{k}\rangle\\
					&+\frac{1}{\eta_{k}}D_{\psi}(x,\bar{x}_{k}). 
				\end{aligned}\label{xk_update}
			\end{eqnarray}
			\STATE Stop if the stopping criterion is reached.
		\end{algorithmic}
	\end{algorithm}	
	
	%%%%%%%%%%%%%%%%
	\begin{remark}
		%From the BPSGE framework in Algorithm \ref{BPSGE}, we have the following results. 
		\begin{itemize}
			\item If $\tilde{\nabla}f(\cdot)=\nabla f(\cdot)$ and $\beta_{k}=0$, the BPSGE algorithm reduces to BPG algorithm \cite{BolteSTV18First}.
			\item If $\tilde{\nabla}f(\cdot)=\nabla f(\cdot)$, the BPSGE algorithm reduces to BPGE, which is a special case of CoCaIn BPG algorithm \cite{MukkamalaOPS20}   and is similar to the BPGE algorithm in  \cite{ZhangBM0C19}.
			\item If  $\beta_{k}=0$, BPSGE reduces to BPSG algorithm \cite{WangH23}. 
			%\red{\item[(iv)] relationship with heavy ball or Nesterov acceleration?? }
		\end{itemize}
	\end{remark}
	
	%%%%%%%%%%%%%
	It is not hard to verify that $\beta_k=0$ always satisfies \eqref{extra_ineq}. 
	%Based on \cite[Lemma 4.1]{MukkamalaOPS20}, we can show that the extrapolation parameter $\beta_{k}$ in \eqref{extra_ineq} exists. 
	%For instance, i
	If $\psi:=\frac{1}{2}\|\cdot\|^{2}$, it follows from \eqref{extra_ineq} that 
	\begin{eqnarray}
		\begin{aligned}
			\beta_{k}^{2}\|x_{k}-x_{k-1}\|^{2}\le& \frac{\delta-\epsilon}{1+\underline{L}\eta_{k-1}}\|x_{k}-x_{k-1}\|^{2}. \end{aligned}\label{example_beta}
	\end{eqnarray}
	Furthermore, it follows from  $\underline{L}=\bar{L}$ and $\delta-\epsilon<1$ that $\beta_{k}< 1/\sqrt{2}$. 
	This insight helps us to choose a proper $\beta_k$ in the numerical experiments.

	%%%%%%%%%%%%%%%%%%%%%%%%%%%%%%%%%%
	%We show the convergence result under the variance reduction assumption similar to  \cite{NguyenLST17, DriggsTLDS2020, WangH23}. 
	An important property in our theoretical analysis is  the variance reduction  given as follows.  
	This definition is similar to that in  \cite{ DriggsTLDS2020, WangH23}.  
	%Now we give the deﬁnition of variance reduced stochastic gradient estimator based on the BPSGE algorithm (Algorithm \ref{BPSGE}) as follows.
	\begin{definition}\label{vr_definition} (Variance reduced stochastic gradient) 
		We say a gradient estimator  $\tilde{\nabla}f(\bar{x}_{k})$ in Algorithm \ref{BPSGE} is variance reduced with constants $V_{1}, V_{2},V_{\Gamma}\ge 0$, and $\tau\in(0,1]$ if it satisfies the following three conditions:
		\begin{itemize}
			\item (MSE (mean squared error) bound): there exist two sequences of random variables $\{\Gamma_{k}\}_{k\ge1}$,  $\{\Upsilon_{k}\}_{k\ge1}$,  such that
			\begin{eqnarray}
				\begin{aligned}
					&\mathbb{E}_{k}[\|\tilde{\nabla}f(\bar{x}_{k})-\nabla f(\bar{x}_{k})\|_{*}^{2}]\\
					\le& \Gamma_{k}+ V_{1}\left(\|x_{k}-x_{k-1}\|^{2}+\|x_{k-1}-x_{k-2}\|^{2}\right),
				\end{aligned}\label{MSE_l22}
			\end{eqnarray}
			and 
			\begin{eqnarray}
				\begin{aligned}
					&\mathbb{E}_{k}[\|\tilde{\nabla}f(\bar{x}_{k})-\nabla f(\bar{x}_{k})\|_{*}]\\
					\le&\Upsilon_{k}+ V_{2}\left(\|x_{k}-x_{k-1}\|+\|x_{k-1}-x_{k-2}\|\right),
				\end{aligned}\label{MSE_l2}
			\end{eqnarray}
			respectively. Here, $\|\cdot\|_{*}$ denotes the conjugate norm \cite{Rockafellar1970} of $\|\cdot\|$.
			\item (Geometric decay): 
			%We can say that 
			The sequence $\{\Gamma_{k}\}_{k\ge1}$ satisfies the following inequality in expectation:
			\begin{eqnarray}
				\begin{aligned}
					&\quad \mathbb{E}_{k}[\Gamma_{k+1}]\le (1-\tau)\Gamma_{k}\\&+V_{\Gamma}\left(\|x_{k}-x_{k-1}\|^{2}+\|x_{k-1}-x_{k-2}\|^{2}\right).\label{Gamma_k1_k}
				\end{aligned}
			\end{eqnarray}
			\item (Convergence of estimator): 
			%For all sequences $\{x_{k}\}_{k=0}^{\infty}$, if 
			If the sequence  $\{x_{k}\}_{k=0}^{\infty}$ satisfies $\lim\limits_{k\rightarrow\infty}\mathbb{E}\|x_{k}-x_{k-1}\|^{2}\rightarrow~0$, then it follows that $\mathbb{E}\Gamma_{k}\rightarrow0$ and $\mathbb{E}\Upsilon_{k}\rightarrow0$.
		\end{itemize}
	\end{definition}

	%%%%%%%%%%%%%%
	\begin{remark}
		\begin{itemize}
			\item In Proposition \ref{vr_gra}, we show   SAGA and SARAH stochastic estimators satisfy Definition \ref{vr_definition}. %\cite{DefazioBL14} \cite{NguyenLST17} 
			\item If $\beta_{k}=0$,   the terms $\|x_{k-1}-x_{k-2}\|^{2}$ and $\|x_{k-1}-x_{k-2}\|$ in \eqref{MSE_l22}, \eqref{MSE_l2} and \eqref{Gamma_k1_k} are not necessary. 
			%See \cite[Definition 3.1]{WangH23} for more details.
			\item We do not require the stochastic gradient to be unbiased or the variance to be bounded. 
		\end{itemize}
	\end{remark}

	%%%%%%%%%%%%%%%%%%%%%%%%%%%%%%%%
	\section{Convergence Analysis}\label{convergence_analysis}
	This section presents a comprehensive analysis of the convergence properties. % exhibited by the BPSGE algorithm. %Firstly, with the general set $C$, we investigate the subsequential convergence of the proposed algorithm and its associated convergence rate in Theorem \ref{subsequence_convergence}. Secondly, we delve into the complexity analysis of the BPSGE algorithm under $C=\mathbb{R}^{d}$. Specifically, the BPSGE algorithm obtains an $\varepsilon$-stationary point of $\Phi$ in expectation in at most $\mathcal{O}(\varepsilon^{-2})$ iterations. For a more thorough understanding of this complexity analysis, please consult Theorem \ref{subgradient_rate}. Finally, under $C=\mathbb{R}^{d}$ and the K{\L} property, we show the global convergence characteristics of the sequence generated by the BPSGE algorithm. For detailed elucidation of these global convergence properties, see Theorem \ref{global_convergence}.  
	%Due to space limit, a
	All proofs are placed in the appendix. 
	
	%%%%%%%%%%%%%%%%%%%%%%%%%%%
	\subsection{Subsequential Convergence Analysis}
	We first show the descent amount of $\Phi(x_{k})$ as follows.
	\begin{lemma}\label{lemma_Phi_kk1}
		Suppose Assumptions \ref{assume_01}-\ref{assume_02} are satisfied  and  $\tilde{\nabla}f(\bar{x}_{k})$  satisfies the variance reduction property defined by Definition \ref{vr_definition}. Let $\{x_{k}\}$ be the sequence generated by Algorithm \ref{BPSGE}. Then the following inequality holds for any $k>0$, 
		\begin{align*}
			&\mathbb{E}_{k}[\Phi(x_{k+1})] +\left(1/\eta_{k}-\alpha-\gamma\right)\mathbb{E}_{k}[D_{\psi}(x_{k},x_{k+1})]\\
			&+\frac{1}{2\tau\gamma}\mathbb{E}_{k}[\Gamma_{k+1}]\\
			\le&\Phi(x_{k})+\frac{1}{2\tau\gamma}\Gamma_{k}+\left(\frac{\delta-\epsilon}{\eta_{k}}+
			\frac{\gamma}{2}\right)D_{\psi}(x_{k-1},x_{k})\\
			&+\frac{\gamma}{2}D_{\psi}(x_{k-2},x_{k-1}).
		\end{align*}
		Here,  $\gamma=\sqrt{2(V_{\Gamma}/\tau+V_{1})}$,  $\alpha$ is the weakly convex parameter in Assumption~\ref{assume_02},  $\delta$ and $\epsilon$ are introduced in \eqref{extra_ineq}, and $V_{1}, V_{2},V_{\Gamma}\ge 0$, $\tau\in(0,1]$ are parameters in Definition~\ref{vr_definition}.
	\end{lemma}
	
	%%%%%%%%%%%%%%%%%%%
	%{\color{blue}Using Lemma \ref{lemma_Phi_kk1}, the following lemma guarantees that $\Psi_{k+1}$ is decreasing in expectation.}
	We introduce a new  Lyapunov function and show it is monotonically decreasing in expectation. 
	\begin{lemma}\label{lyapunov_descent}
		Suppose the same conditions with Lemma \ref{lemma_Phi_kk1} hold, and  the stepsize satisfies 
		\begin{eqnarray}
			\eta_{k}\le \min\left\{\eta_{k-1}, \frac{1}{\bar{L}}, \frac{1-\delta}{\alpha+2\gamma}\right\}, \,\, \forall\, k>0. \label{stepsize_set}
		\end{eqnarray}
		Let $\{x_{k}\}_{k\in\mathbb{N}}$ be a sequence generated by BPSGE (Algorithm \ref{BPSGE}) and define the following Lyapunov sequence 
		\begin{eqnarray}
			\begin{aligned}
				\Psi_{k+1}:= &\eta_{k}(\Phi(x_{k+1}) -\mathcal{V}(\Phi))  +t_{k}D_{\psi}(x_{k},x_{k+1})\\
				&+\eta_{k}\left(\frac{\gamma}{2}+\frac{\epsilon}{3\eta_{k}}\right)D_{\psi}(x_{k-1},x_{k})+\frac{\eta_{k}\Gamma_{k+1}}{2\tau\gamma}, \label{lyapunov_function}
			\end{aligned}
		\end{eqnarray}
		where $t_{k}=1-\eta_{k}\alpha-\eta_{k}\gamma-\epsilon/3$. Then, for all $k\in\mathbb{N}$, we have
		\begin{eqnarray}
			\begin{aligned} %\Psi_{k}-\mathbb{E}_{k}[\Psi_{k+1}]			\ge
				\mathbb{E}_{k}[\Psi_{k+1}]&\le\Psi_{k}-\frac{\epsilon}{3}(\mathbb{E}_{k}[D_{\psi}(x_{k},x_{k+1})]\\
				&+D_{\psi}(x_{k-1},x_{k})+D_{\psi}(x_{k-2},x_{k-1})).
			\end{aligned}\label{decent_inequality_01}
		\end{eqnarray} 
	\end{lemma}
	
	%%%%%%%%%%%%%%%%%%%%%%%%%%
	Based on Lemma \ref{lyapunov_descent}, we get the subsequential convergence of BPSGE.
	\begin{theorem}\label{subsequence_convergence}
		Let $\{x_{k}\}_{k\in\mathbb{N}}$ be a sequence generated by BPSGE algorithm. Then, the following statements hold.
		\begin{itemize}
			\item The sequence $\{\mathbb{E}[\Psi_{k}]\}_{k\in\mathbb{N}}$ is nonincreasing.
			\item $\sum\limits_{k=1}^{+\infty}\mathbb{E}[D_{\psi}(x_{k-1},x_{k})]<+\infty$,
			and    the sequence $\{\mathbb{E}[ D_{\psi}(x_{k-1}, x_{k}) ]\}$  %$\{\mathbb{E}[ D_{\psi}(x_{k}, x_{k+1}) ]\}$ 
			converges to zero. 
			\item $\min\limits_{1\le k\le K}\mathbb{E}[D_{\psi}(x_{k-1},x_{k})]\le 3\Psi_{1}/(\epsilon K)$.
		\end{itemize}
	\end{theorem}
	
	%%%%%%%%%%%%%%%%%%%%
	\subsection{Global Convergence Analysis}
	Now we show the convergence of  the whole sequence to a  stochastic stationary point under more conditions.
	Consider the case  $C=\mathbb{R}^{d}$. 
	%and analyze the convergence rate of the expected squared distance of the subgradient and global convergence of BPSGE. 
	We require the following additional assumptions. 
	\begin{assumption}\label{assume_03}
		\begin{itemize}
			\item $\nabla f_{i}(x)$ is Lipschitz continuous with constant $M_{1} > 0$ on any bounded subset of $\mathbb{R}^{d}$. 
			\item   $\nabla \psi $ is Lipschitz continuous with  constant $M_{2} > 0$ on any bounded subset of $\mathbb{R}^{d}$.
		\end{itemize}
	\end{assumption}
	%%%%%%%%%%%%%%%%%%%%%%%
	From Assumption \ref{assume_03}\,(1), the function $\nabla f(x)$ is also Lipschitz continuous with constant $M_{1} > 0$ on any bounded subset of $\mathbb{R}^{d}$. 
	%Under Definition \ref{vr_definition}, and the definition of  
	Furthermore, we show SAGA \cite{DefazioBL14} and SARAH \cite{NguyenLST17} satisfy  the following proposition.
	\begin{proposition}\label{vr_gra}
		Assume Assumption \ref{assume_03}\,(1)  holds and the mini-batch $B_{k}$ is uniform randomly chosen from   $\{ 1, \dots , n \}$ with cardinality $b$, i.e., $b=|B_{k}|$. 
		\begin{itemize}
			\item The SAGA gradient estimator  \cite{DefazioBL14} 
			%\begin{eqnarray*}
			\begin{align*}
				\tilde{\nabla}^{SAGA}f(\bar{x}_{k})=&\frac{1}{b}(\sum_{j\in B_{k}}\nabla f_{j}(\bar{x}_{k})-g_{k}^{j})+\frac{1}{n}\sum_{i=1}^{n}g_{k}^{i},\\
				g_{k+1}^{i}=&
				\begin{cases}
					\nabla f_{i}(\bar{x}_{k}), &\text{if }  i\in B_{k}, \\
					g_{k}^{i},&\text{otherwise.}
				\end{cases}	
			\end{align*}%\label{SAGA_gra}
			%\end{eqnarray*}
			is variance reduced with
			\begin{align*}
				\Gamma_{k+1}:=&\frac{1}{bn}\sum_{i=1}^{n}\|\nabla f_{i}(\bar{x}_{k})-\nabla f_{i}(z_{k}^{i})\|_{*}^{2},\\
				\Upsilon_{k+1}:=&\frac{1}{\sqrt{bn}}\sum_{i=1}^{n}\|\nabla f_{i}(\bar{x}_{k})-\nabla f_{i}(z_{k}^{i})\|_{*},
			\end{align*}
			where $z_{k}^{i}=x_{k}$ if $i\in B_{k}$ and $z_{k}^{i}=z_{k-1}^{i}$ otherwise. The constants $\tau=\frac{b}{2n}$, $V_{\Gamma}=\frac{2b+4n}{b^{2}}M_{1}^{2}$, $V_{1}=V_{2}=M_{1}$ $V_{1}=M_{1}^2, V_{2}=M_{1}$.
			\item The SARAH gradient estimator \cite{NguyenLST17}  
			\begin{eqnarray*}
				&&\tilde{\nabla}^{SARAH}f(\bar{x}_{k})\\
				&=&
				\left\{\begin{array}{ll}
					\nabla f(\bar{x}_{k}),\,\,&  \mbox{w.p.}\,\,\frac{1}{p},\\
					\frac{1}{b}(\underset{j\in B_{k}}{\sum}\nabla f_{j}(\bar{x}_{k})-\nabla f_{j}(\bar{x}_{k-1}))\\
					\quad+\tilde{\nabla}^{SARAH}f(\bar{x}_{k-1}),&\mbox{otherwise.}
				\end{array}\right.
			\end{eqnarray*}
			is variance reduced with
			\begin{align*}
				\Gamma_{k+1}=&\|\tilde{\nabla}^{SARAH}f(\bar{x}_{k})-\nabla f(\bar{x}_{k})\|_{*}^{2},\\
				\Upsilon_{k+1}=&\|\tilde{\nabla}^{SARAH}f(\bar{x}_{k})-\nabla f(\bar{x}_{k})\|_{*},
			\end{align*}
			and constants $\tau=\frac{1}{p}, V_{1}=V_{\Gamma}=2M_{1}^{2}, V_{2}=2M_{1}$. Here ``w.p. $\frac{1}{p}$'' means with probability $\frac{1}{p} \in (0, 1]$.  
		\end{itemize}
	\end{proposition}
	%%%%%%%%%%%%%%%%%%%%%%
	\begin{corollary}\label{lemma_Phi_kk1_Rd}
		If $C=\mathbb{R}^{d}$,  from Lemma \ref{lyapunov_descent}, it shows that
		\begin{align*}
			\mathbb{E}_{k}[\Psi_{k+1}]&\le \Psi_{k}-\frac{\epsilon}{6}(\mathbb{E}_{k}\|x_{k+1}-x_{k}\|^{2}\\
			&+\|x_{k}-x_{k-1}\|^{2}+\|x_{k-1}-x_{k-2}\|^{2}).
		\end{align*}
	\end{corollary}
	Now we prove the following bound for   $\partial\Phi(x_{k})$.
	\begin{lemma} \label{subgradient_bound}
		Suppose that Assumptions \ref{assume_01}-\ref{assume_03} hold and the stepsize $\eta_{k}$ satisfies $0<\eta\le \eta_{k}$\footnote{$\eta$ is the lower bound of $\eta_{k}$.} and \eqref{stepsize_set}. Let $\{x_{k}\}_{k\in\mathbb{N}}$ be a bounded sequence generated by the BPSGE algorithm. Define
		\[
		w_{k+1}:=\nabla f(x_{k+1})-\tilde{\nabla}f(\bar{x}_{k})+\frac{1}{\eta_{k}}(\nabla \psi(\bar{x}_{k})-\nabla \psi(x_{k+1})).
		\]
		Then, we have $w_{k+1}\in\partial \Phi(x_{k+1})$ and 
		\begin{align*}
			&\mathbb{E}_{k}\|w_{k+1}\|
			\le\rho(\mathbb{E}_{k}\|x_{k+1}-x_{k}\|+\|x_{k}-x_{k-1}\|\\
			&+\|x_{k-1}-x_{k-2}\|) +\Upsilon_{k},
		\end{align*}
		where $\rho=\max\{M_{1}+\frac{M_{2}}{\eta}, V_{2}+\beta_{k}M_{1}+\frac{\beta_{k}M_{2}}{\eta}, V_{2}\}$.
	\end{lemma}
	
	%%%%%%%%%%%%%%%
	Similarly,  we show the   bound for   $\mathbb{E}[\mathrm{dist}(0,\partial\Phi(x_{k+1}))^{2}]$.
	%Similarly, we get the following results.
	\begin{lemma}\label{lemma_dist2}
		Under the same conditions as in Lemma \ref{subgradient_bound}, there exists a constant $\bar{\rho}>0$ such that
		\begin{align*}
			&\mathbb{E}[\mathrm{dist}(0,\partial\Phi(x_{k+1}))^{2}]
			\le\bar{\rho}\mathbb{E}[\|x_{k+1}-x_{k}\|^{2}\\
			&+\|x_{k}-x_{k-1}\|^{2}+\|x_{k-1}-x_{k-2}\|^{2}] +3\mathbb{E}\Gamma_{k}.
		\end{align*}
	\end{lemma}
	
	%%%%%%%%%%%
	Using Lemma \ref{lemma_dist2}, we show the $\mathcal{O}(1/\epsilon^{2})$ complexity  in expectation  to obtain an $\epsilon$-stationary point. 
	%show the $O(\frac1k)$ convergence rate of the expected squared distance of the subgradient to zero.
	\begin{theorem}\label{subgradient_rate}
		Assume that Assumptions \ref{assume_01}-\ref{assume_03} hold, and the stepsize satisfies $0<\eta\le\eta_{k}$ and \eqref{stepsize_set}.  Let $\{x_{k}\}_{k\in\mathbb{N}}$ be a bounded sequence generated by the BPSGE algorithm. Then there exists  some $0<\sigma<\epsilon/6$ such that
		\begin{align*}
			\mathbb{E}[\mathrm{dist}(0,\partial\Phi(x_{\hat{k}}))^{2}]\le& \frac{\bar{\rho}}{(\epsilon/6-\sigma)K}(\mathbb{E}\Psi_{1}+\frac{\epsilon/2-3\sigma}{\tau\bar{\rho}}\mathbb{E}\Gamma_{1})\\
			=&\mathcal{O}(1/K),
		\end{align*}
		where $\hat{k}$ is drawn from $\{2, \dots, K+1\}$. In other words, it takes at most $\mathcal{O}(1/\epsilon^{2})$ iterations  in expectation  to obtain an $\epsilon$-stationary point (Definition \ref{stationary-point}) of $\Phi$. 
	\end{theorem}
	
	We define the set of limit points of $\{x_{k}\}_{k=0}^{\infty}$ as
	\begin{align*}
		\omega(x_{0}):=\{&x: \exists \text{ an increasing sequence of integers } \{k_{l}\}_{l\in\mathbb{N}}\\
		&\text{ s.t. } x_{k_{l}}\rightarrow x\text{ as } l\rightarrow \infty\}. 
	\end{align*}
	Now we get the properties of limit points of $\{x_{k}\}$ as follows. 
	\begin{lemma}\label{statements_lemma}
		Suppose that Assumptions \ref{assume_01} to \ref{assume_03} hold,   the step $\eta_{k}$ satisfies $0<\eta\le \eta_{k}$ and \eqref{stepsize_set}.  Then the following statements hold. 
		\begin{itemize}
			\item $\sum_{k=0}^{\infty}\|x_{k+1}-x_{k}\|^{2}<+\infty$ a.s., and $\lim_{k\rightarrow+\infty}\|x_{k+1}-x_{k}\|\rightarrow0$ a.s.
			\item $\mathbb{E}[\Phi(x_{k})]\rightarrow\Phi_{*}$, where $\Phi_{*}\in[\bar{\Phi},+\infty)$  with $\bar{\Phi}:=\inf_{x} \Phi(x)$, and $\mathbb{E}\Phi(x_{*})=\Phi_{*}$ for all $x_{*}\in\omega(x_{0})$.
			\item $\mathbb{E}[\text{dist}(0,\partial\Phi(x_{k}))]\rightarrow0$. Moreover, the set $\omega(x_{0})$ is nonempty, and $\mathbb{E}[\text{dist}(0,\partial\Phi(x_{*}))]=0 \text{ for all } x_{*}\in\omega(x_{0})$.
			\item $\text{dist}(x_{k},\omega(x_{0}))\rightarrow0$ a.s., and $\omega(x_{0})$ is a.s. compact and connected.
		\end{itemize}
	\end{lemma}	
	
	%Based on Lemma 6 in the Appendix,  Corollary \ref{lemma_Phi_kk1_Rd}, and Lemma \ref{subgradient_bound},  we give the global convergence result of the BPSGE algorithm.
	We further show the whole sequence convergence under the K{\L} property. 
	\begin{theorem}\label{global_convergence}
		Suppose that Assumptions \ref{assume_01}-\ref{assume_03} hold, the step $\eta_{k}$ satisfies $0<\eta\le \eta_{k}$ and \eqref{stepsize_set}. Let $\{x_{k}\}_{k\in\mathbb{N}}$ be a bounded sequence generated by the BPSGE algorithm. If the optimization function $\Phi(x)$ is a semialgebraic function that satisﬁes the K{\L} property with exponent $\theta \in [0, 1)$ (see Lemma 6 in Appendix), then either the point $x_{k}$ is a critical point after a ﬁnite number of iterations or the sequence $\{ x_{k}\}_{k \in\mathbb{N}}$ almost surely satisﬁes the ﬁnite length property in expectation, namely,
		\[
		\sum_{k=0}^{+\infty}\mathbb{E}\|x_{k+1}-x_{k}\|<+\infty.
		\]
	\end{theorem}
	%%%%%%%%%%%%%%%%%%%%%%%%%%%%%%%%%%%%%%%%%%%%%%%%%%%%%%%

	\section{Numerical Experiments}\label{numercial_experiments}
	%%%%%%%%%%%%%%%%%%%%%%%%%
	In this section, we present our numerical study on the practical performance of the proposed BPSGE algorithm with three diﬀerent stochastic gradient estimators.  For all numerical experiments,  $C=\mathbb{R}^{d}$.  The experiments are implemented in Matlab 2020b and conducted on a computer with AMD Ryzen 5 5600x 6-Core 3.7 GHz and 48GB memory.
	
	For the stochastic algorithms, we repeat all numerical experiments $10$ times and report their average performance. All  the initial points are generated by the uniform distribution between $0$ and $0.1$ and are the same for all algorithms. Inspired by \eqref{example_beta}, we set $\beta^{k}=0.6\frac{k-1}{k+2}$  for simplicity\footnote{The inequality \eqref{extra_ineq} in Algorithm \ref{BPSGE} is required for our convergence analysis. It is time-consuming to check this inequality in numerical experiments. Due to this reason, we directly set $\beta^{k}=0.6\frac{k-1}{k+2}$. Our numerical experiments show that BPSGE always converges with this $\beta^k$.} in the BPGE and BPSGE-SGD/SAGA/SARAH algorithms. The stepsize is set as $\eta_{k}=\min(\eta_{k-1}, L_{k}^{-1})$, where $L_{k}$ is the approximated Lipschitz constant 
	estimated by the power method to $\bar{V}_{k}$ and $\bar{U}_{k}$ for $U$- and $V$-update,  respectively. 
	This choice of stepsize is  the  same  for all  compared algorithms. 
	%The source codes are available at our Github repository\footnote{\url{https://github.com/nothing2wang/BPSGE-Algorithm}}. 

	%%%%%%%%%%%%%%%%%%%%%%%%%%%%%%%%%%%%%%%%%%%%%%%%
	%%%%%%%%%%%%%%%%%%%%
	\begin{table*}[!ht]
		\centering
		%\setlength{\tabcolsep}{1pt}
		%\fontsize{8}{11}\selectfont
		\begin{adjustbox}{max width=\textwidth}
			\begin{tabular}{c|c|c|c|c|c|c|c|c}\hline
				%\multirow{3}{*}{Dataset}&\multicolumn{8}{c|}{Accuracy ($\%$)}\\\cline{2-9}
				%& BPG &BPSG-SGD & BPSG-SARAH  & BPSG-SAGA& BPGE& BPSGE-SGD&BPSGE-SARAH & BPSGE-SAGA \\\hline
				\multirow{2}{*}{Dataset}&\multirow{2}{*}{BPG} & \multicolumn{3}{c|}{BPSG} & \multirow{2}{*}{BPGE} & \multicolumn{3}{c}{BPSGE} \\ \cline{3-5}\cline{7-9}
				& & -SGD & -SARAH & -SAGA &  & -SGD & -SARAH & -SAGA \\\hline
				\emph{COIL20} & 75.56 &79.30& 85.47 & 85.83& 85.33& 87.48&89.37&\textbf{89.97}\\
				\emph{PIE} & 77.21 &84.85& 86.03 & 86.33& 84.27& 87.79&88.45&\textbf{88.83}\\
				\emph{COIL100} & 73.06 &80.70& 82.23 & 82.46& 76.25& 82.24&84.02&\textbf{84.29}\\
				\emph{TDT2} & 68.08 &82.01& 85.04 & 85.22& 77.92& 85.21&87.42&\textbf{87.54}\\\hline
			\end{tabular}
		\end{adjustbox}
		\caption{Comparison of clustering accuracy ($\%$)   on four datasets by graph regularized NMF.}
		\label{GNMG_clustering}
	\end{table*}
	
	\subsection{Graph Regularized NMF for Clustering}
	Previous studies \cite{Cai11GNMF, ShahnazBPP06, AhmedHAD21}  show that NMF is  powerful for clustering problems, especially in document clustering and image clustering tasks. 
	%For example, a face image can be thought of as a combination of nose, mouth, eyes, etc. It is also reasonable to require the combination coefficients to be non-negative. This is the main motivation for applying NMF to image clustering.  
	We consider the   graph regularized NMF problem for clustering proposed in \cite{Cai11GNMF} which is defined as 
	\begin{eqnarray}
		\min_{U\in\mathbb{R}^{m\times r}_{+},V\in\mathbb{R}^{r\times d}_{+}} \frac{1}{2}\|M-UV\|_{F}^{2}+\frac{\mu_{0}}{2}\text{Tr}(U^{T}LU), \label{GNMF}
	\end{eqnarray}
	where   $\|\cdot\|_{F}$ denotes the Frobenius norm, $\text{Tr}(\cdot)$ denotes the trace of a matrix, $L$ denotes the graph Laplace matrix,  and $\mu_{0}$ is a positive parameters. 
	%The graph Laplacian matrix $L$ satisfies $L= D-S$, where $D$ is a diagonal   degree matrix \cite{Cai11GNMF, CaiHHZ06} and $S$ is the adjacency matrix. 
	This model can help to distinguish anomalies from normal observation \cite{AhmedHAD21}. 
	We first define the following  kernel-generating distance, i.e.,
	\begin{eqnarray}
		\begin{aligned}
			\psi_{1}(U,V):=&\left(\|U\|_{F}^{2}/2+\|V\|_{F}^{2}/2\right)^{2},\\ \psi_{2}(U,V):=&\|U\|_{F}^{2}/2+\|V\|_{F}^{2}/2, 
		\end{aligned}\label{kernels}
	\end{eqnarray}
	which is designed to allow for closed-form update in the BPSGE algorithm. 
	Let 
	\begin{eqnarray}
		\begin{aligned}
			f(U,V):=&\frac{1}{2}\|M-UV\|_{F}^{2}+\frac{\mu_{0}}{2}\text{Tr}(U^{T}LU),\\
			h(U,V):=&I_{U\ge 0}+I_{V\ge 0},\\
			\psi(U,V):=&3\psi_{1}(U,V)+(\|M\|_{F}+\mu_{0}\|L\|_{F})\psi_{2}(U,V),
		\end{aligned}\label{fh_def_01}
	\end{eqnarray}
	where $I_{U\ge0}$ is the indicator function. 
	Now we give the closed form of $(U_{k+1},V_{k+1})$ as follows.
	\begin{proposition}\label{Prop_GNMF}
		With the above defined $f,h,\psi$ in \eqref{fh_def_01}, the update steps \eqref{xk_update} in each iteration are given by 
		\[
		U_{k+1}=t\Pi_{+}(-P_{k}),\quad V_{k+1}=t\Pi_{+}(-Q_{k}),
		\] 
		where 
		\begin{align}\label{eqn:pkqk}
			\begin{aligned}
				P_{k}=&\eta_{k}\tilde{\nabla}_{U}f(\bar{U}_{k},\bar{V}_{k})-\nabla_{U}\psi(\bar{U}_{k},\bar{V}_{k}),\\
				Q_{k}=&\eta_{k}\tilde{\nabla}_{V}f(\bar{U}_{k},\bar{V}_{k})-\nabla_{V}\psi(\bar{U}_{k},\bar{V}_{k}).
			\end{aligned}
		\end{align} 
		Here, $\Pi_{+}(\cdot)$ is the projection onto the nonnegative orthant,  and $t\ge0$   satisfies
		\begin{align*}
			&3\left(\|\Pi_{+}(-P_{k})\|_{F}^{2}+\|\Pi_{+}(-Q_{k})\|_{F}^{2}\right)t^{3}\\
			&+(\|M\|_{F}+\mu_{0}\|L\|_{F})t-1=0.
		\end{align*}
	\end{proposition}
	
	We use four datasets \emph{COIL20}, \emph{PIE}, \emph{COIL100}, and \emph{TDT2} \footnote{\url{http://www.cad.zju.edu.cn/home/dengcai/Data/data.html}}   to illustrate the numerical performance of the proposed algorithm. In this numerical experiment, we let $\mu_{0}=100$ , and $r=20$ for \emph{COIL20}, $r=68$  for \emph{PIE}, $r=100$ for \emph{COIL100}, and $r=30$ for \emph{TDT2}, respectively. 
	%In this study, we conduct a series of experiments and compare the performance of the four datasets over 50 epochs. To ensure controlled and efficient computations, we employ a minibatch subsampling ratio set at $5\%$.
	Here we conduct experiments  with 50 epochs using the minibatch subsampling ratio  $5\%$.
	
	% \begin{table}[!ht]
	% 	\centering
	% 	%\setlength{\tabcolsep}{3pt}
	% 	\fontsize{9}{10}\selectfont
	% 	\begin{tabular}{c|c|c|c}\hline
	% 		Dataset& Size &Dimensionality     & Number of classes \\\hline
	% 		\emph{COIL20} & 1440 & 1024 &20\\
	% 		\emph{PIE} & 2856 &1024 & 68\\
	% 		\emph{COIL100} & 7200 & 1024 &100\\
	% 		\emph{TDT2} & 9394 & 36771 & 30\\\hline
	% 	\end{tabular}
	%  \caption{Statistics of the four datasets in graph regularized NMF.}
	% 	\label{GNMG_datasets}
	% \end{table}
	
	After solving \eqref{GNMF} by Proposition~\ref{Prop_GNMF}, we compute the clustering label by implementing the $K$-means method \cite{litekmeans} on $U$. 
	%The clustering accuracy is the ratio of correctly predicted class labels. 
	%comparing the obtained label and the true label \cite{CaiHH05}.  
	All results are presented in Table \ref{GNMG_clustering}. 
	%of each sample with the label provided by the data set. 
	%The accuracy (AC) is used to measure the clustering performance. 
	%Please see \cite{CaiHH05} for more details. 
	%Table \ref{GNMG_clustering} shows the clustering results on four datasets for all compared algorithms, respectively. 
	From this table, it shows that the extrapolation technique can improve the numerical performance. In addition, the stochastic algorithms can get better numerical results than their deterministic versions. Furthermore, the variance reduction stochastic gradient estimators can get the best performance in the stochastic framework.  
	
	\subsection{MF with Weakly-convex Regularization}\label{Sparse_NMF}
	%Given a matrix $M\in\mathbb{R}^{m\times d}$, to obtain the factors $U\in\mathbb{R}^{m\times r}$ with a weakly-convex regularization \cite{YinLHX15, MaLH17} and $V\in\mathbb{R}^{r\times d}$  (where $r<\min\{m, d\}$), we solve 
	
	Consider the following  optimization problem with a weakly-convex regularization \cite{YinLHX15, MaLH17} 
	\begin{eqnarray}
		\min_{U,V}\frac{1}{2}\|M-UV\|_{F}^{2}+\lambda_{1}\|U\|_{1}-\frac{\lambda_{2}}{2}\|U\|_{F}^{2},\label{WCMF}
	\end{eqnarray}
	where $U\in\mathbb{R}^{m\times r},V\in\mathbb{R}^{r\times d}$, $\|U\|_{1}:=\sum_{i,j}|U_{i,j}|$ and $\lambda_{1},\lambda_{2}$ are two positive parameters. The term $\lambda_{1}\|U\|_{1}-\frac{\lambda_{2}}{2}\|U\|_{F}^{2}$ is a $\lambda_{2}$-weakly convex function. 
	
	%The  model \eqref{WCMF} can be covered by the optimization  problem \eqref{target_function} with $f(U,V):=\frac{1}{2}\|M-UV\|_{F}^{2}$ and $h(U,V):=\lambda_{1}\|U\|_{1}-\frac{\lambda_{2}}{2}\|U\|_{F}^{2}$.  
	Now we give the closed-form of $(U_{k+1}, V_{k+1})$ for the problem \eqref{WCMF}  in the following proposition.
	
	\begin{proposition}\label{prop_WCMF}
		If $f(U,V):=\frac{1}{2}\|M-UV\|_{F}^{2}$, $h(U,V):=\lambda_{1}\|U\|_{1}-\frac{\lambda_{2}}{2}\|U\|_{F}^{2}$, $\psi(U,V):=3\psi_{1}(U,V)+\|M\|_{F}\psi_{2}(U,V)+\frac{\eta\lambda_{2}}{2}\|U\|_{F}^{2}$ (where $\eta=\eta_{k}$ in the $k$-th iteration), we have  the update steps  for solving \eqref{xk_update}  in each iteration  are 
		\[
		U_{k+1}=t\mathcal{S}_{\lambda_{1}\eta_{k}}(-P_{k}),\quad V_{k+1}=-tQ_{k},
		\] 
		where $P_k$ and $Q_k$ are defined by \eqref{eqn:pkqk}, $t\ge0$ satisfies 
		\[
		3(\|\mathcal{S}_{\lambda_{1}\eta_{k}}(-P_{k})\|_{F}^{2}+\|-Q_{k}\|_{F}^{2}) t^{3}+\|M\|_{F}t-1=0,
		\]
		and $\mathcal{S}_{\lambda_{1}\eta_{k}}(\cdot)$ is the soft-thresholding operator\footnote{For any $y\in\mathbb{R}^{d}$, $S_{\tau}(y) =\underset{x\in\mathbb{R}^{d}}{\arg\min}\{\tau\|x\|_{1}+\frac{1}{2}\|x-y\|^{2}\}=\max\{|y|-\tau,0\}\text{sign}(y)$.} \cite{Donoho95}.
	\end{proposition}

	We use two real datasets\footnote{\url{http://www.cad.zju.edu.cn/home/dengcai/Data/FaceData.html}}, i.e., \emph{ORL} with the size of $4096 \times 400$ and \emph{Yale-B} with the size of $2414 \times 1024$, to illustrate the numerical performance.   We let $\lambda_{1}=0.05$ and $\lambda_{2}=0.02$ for all numerical experiments, and let $r=25$ for \emph{ORL} dataset and $r=49$ for \emph{Yale-B} dataset, respectively. 
	%Here we conduct experiments to compare BPG, BPSG-SGD/SAGA/SARAH, BPGE, and BPSGE-SGD/SAGA/SARAH for \emph{ORL} and \emph{Yale-B}  datasets
	We conduct experiments with 200 epochs using the minibatch subsampling ratio  $5\%$. 
	%The numerical results in Fig.~\ref{WCMF_epochs}  show  the mean of log of the objective function $\Phi(U_{k}, V_{k})$ versus the number of epochs   for \emph{ORL} and \emph{Yale-B} datasets, respectively. 
	\begin{figure}[!t]
		\centering
		\includegraphics[width=0.45\textwidth]{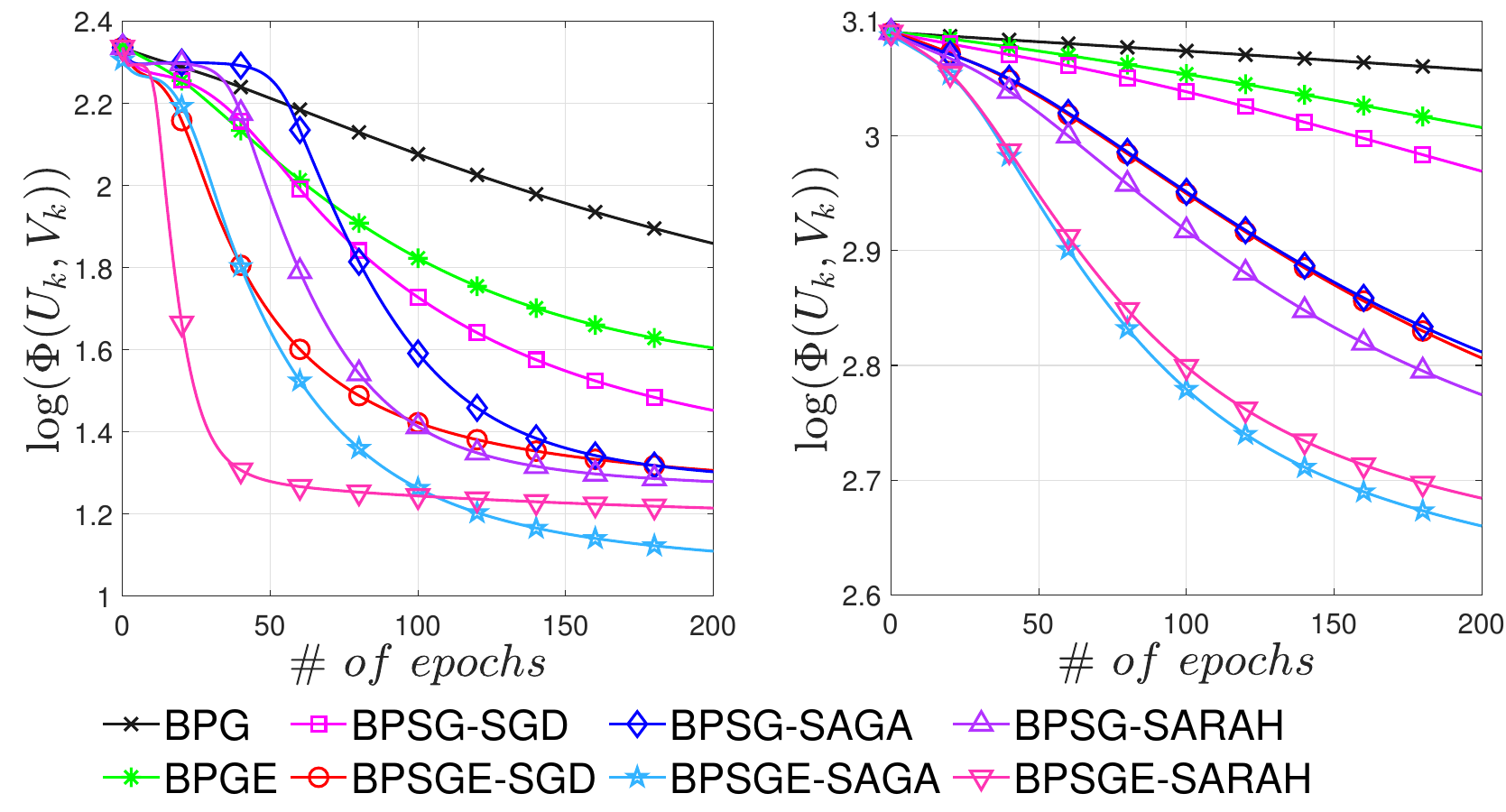}
		\caption{Numerical experiment results on \emph{ORL} and \emph{Yale-B} datasets for   problem \eqref{WCMF}. Left: \emph{ORL} with $r=25$. Right: \emph{Yale-B} with $r=49$.
			%Numerical experiment results on \emph{ORL} and \emph{Yale-B} datasets for optimization problem \eqref{WCMF}: log of objection value $\Phi(U_{k}, V_{k})$ versus epochs by different algorithms.
		}
		\label{WCMF_epochs}
	\end{figure}
	%The numerical results shown in 
	% Fig.~\ref{WCMF_epochs} indicates that the extrapolation-based algorithms perform better than those without extrapolation. Additionally, the stochastic gradient methods produce lower objective function values compared to their deterministic counterparts. Furthermore, the variance reduced stochastic gradient algorithms (with SAGA/SARAH) show better performance compared to SGD-based methods.
	
	%The numerical results shown in 
	Fig.~\ref{WCMF_epochs} indicates that the extrapolation-based algorithms perform better than those without extrapolation. %Additionally, the stochastic gradient methods produce lower objective function values compared to their deterministic counterparts. 
	Furthermore, the variance reduced stochastic gradient algorithms (with SAGA/SARAH) show better performance compared to SGD-based methods.

	\subsection{NMF with Nonconvex Sparsity Constraints}
	We continue to study the NMF problem with nonconvex sparsity constraints given by 
	\begin{eqnarray}
		\min_{U,V} \left\{\frac{1}{2}\|M-UV\|_{F}^{2}: \|U_{:,i}\|_{0}\le s_{1}, \|V_{j,:}\|_{0}\le s_{2}\right\}, \label{SSNMF}
	\end{eqnarray}
	where $U\in\mathbb{R}^{m\times r}_{+},V\in\mathbb{R}^{r\times d}_{+}$, $r>0$, $i,j\in\{1,2\dots,r\}$, $U_{:,i}$ denotes the $i$-th column of $U$ and $\|U_{:,i}\|_{0}$ denotes the number of non-zero entries of $U_{:,i}$. Similarly,  $V_{j,:}$ is the $j$-th row of $V$. Now, we denote 
	\begin{align*}
		f(U,V):=&\frac{1}{2}\|M-UV\|_{F}^{2},\\
		h(U,V):=&I_{U\ge0}+I_{\|U_{:,1}\|_{0}\le s_{1}}+\dots+I_{\|U_{:,r}\|_{0}\le s_{1}}\\&+I_{V\ge 0}+I_{\|V_{1,:}\|_{0}\le s_{2}}+\dots+I_{\|V_{r,:}\|_{0}\le s_{2}},\\
		\psi(U,V): =&3\psi_{1}(U,V)+\|M\|_{F}\psi_{2}(U,V).
	\end{align*}
	Here, $\psi_1$ and $\psi_2$ are given by \eqref{kernels}. 
	The closed-form of $(U_{k+1},V_{k+1})$   is the same as that of   \cite[Proposition D.8]{MukkamalaO19} and is shown as follows. 
	
	%We present this closed form in the following proposition and omit the proof.
	\begin{proposition}\label{Prop_SSNMF}
		Given the optimization problem \eqref{SSNMF} with the above defined $f(\cdot)$, $h(\cdot)$ and $\psi(\cdot)$, the update  \eqref{xk_update} in the BPSGE algorithm   are given by 
		\[
		U_{k+1}=t\mathcal{H}_{s_{1}}(\Pi_{+}(-P_{k})), \quad V_{k+1}=t\mathcal{H}_{s_{2}}(\Pi_{+}(-Q_{k})),
		\]
		where $P_k$ and $Q_k$ are defined by \eqref{eqn:pkqk}, $t\ge 0$ and satisfies 
		\[\begin{aligned}
			3(\|\mathcal{H}_{s_{1}}(\Pi_{+}(-P_{k}))\|_{F}^{2}&+\|\mathcal{H}_{s_{2}}(\Pi_{+}(-Q_{k}))\|_{F}^{2})t^{3} \\&+\|M\|_{F}t-1=0.
		\end{aligned}\]
		Here, $\Pi_{+}(\cdot)$ is the projection on the nonnegative space, $\mathcal{H}_{s}(\cdot)$ is the hard-thresholding operator\footnote{Let $y\in\mathbb{R}^{d}$ and   $\vert y_{1}\vert \ge \vert y_{2}\vert \ge \dots  \ge \vert y_{d}\vert$. Then the hard-thresholding operator  is given by
			\[
			\mathcal{H}_{s}(y)=\arg\min_{x\in\mathbb{R}^{d}} \{\|x-y\|^{2}:\|x\|_{0}\le s\}=\begin{cases}
				y_{i},  & i\le s,\\
				0, & \text{otherwise.}
			\end{cases}
			\]
			where $s > 0$ and the operations are applied element-wise.} \cite{LussT13}.
	\end{proposition}
	
	\begin{figure}[!ht]
		\centering
		\begin{tabular}{cc}
			\includegraphics[width=0.45\textwidth]{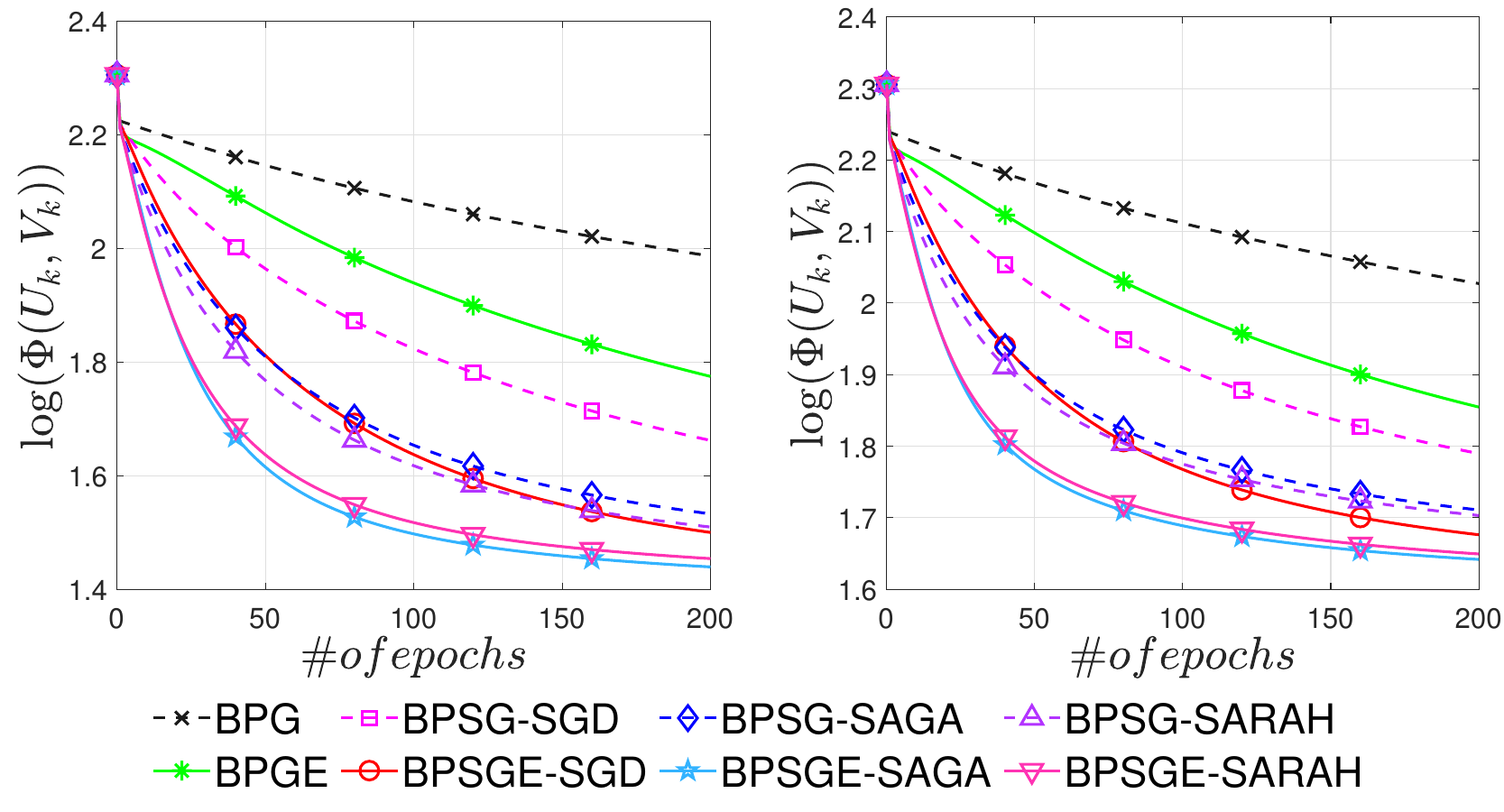}\\
			(a) Left: $s_{1}=\frac{m}{2}$ and $s_{2}=\frac{d}{3}$. Right:  $s_{1}=\frac{m}{2}$ and $s_{2}=\frac{d}{5}$.\\
			\includegraphics[width=0.45\textwidth]{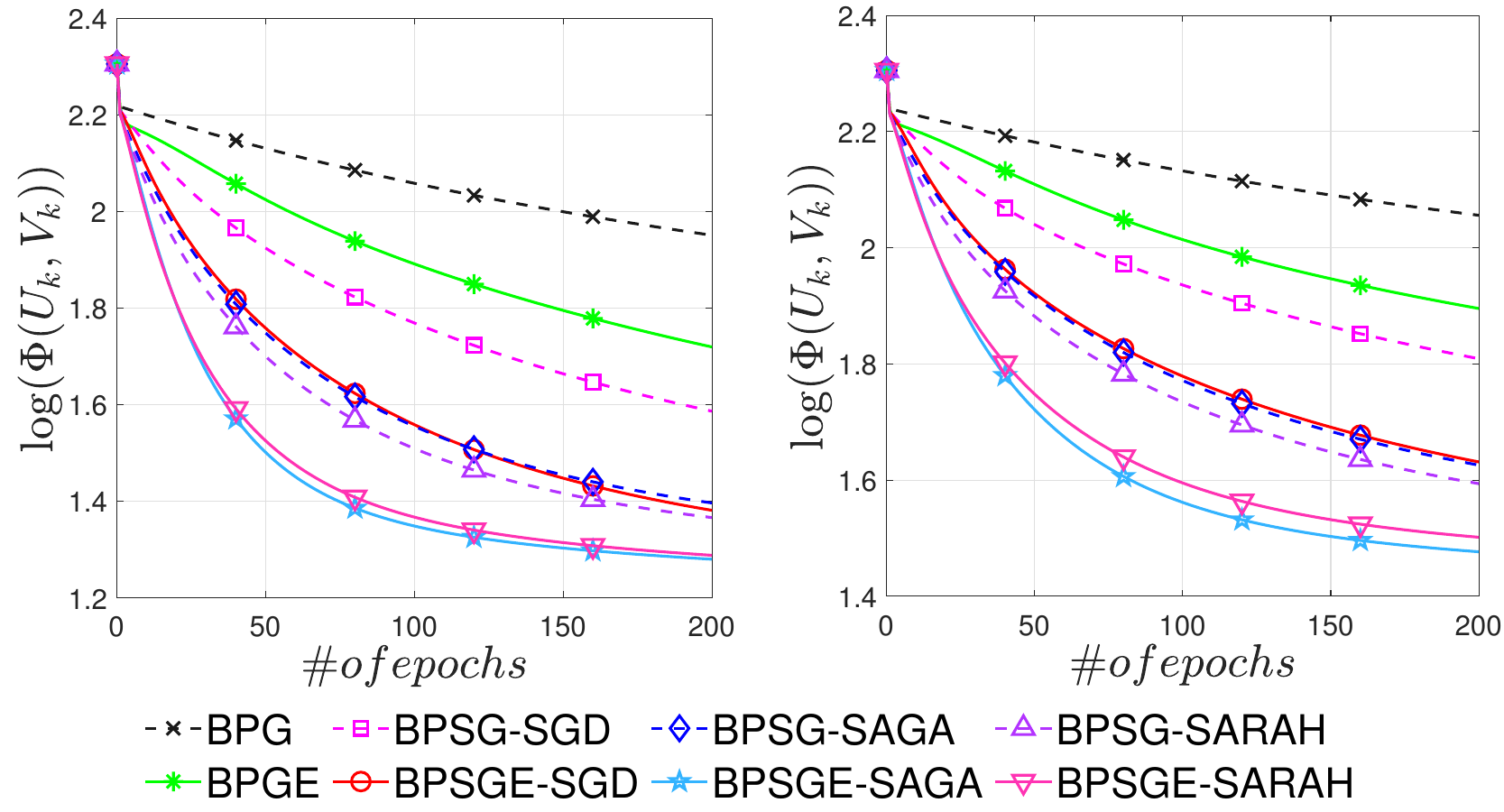}\\
			(b) Left: $s_{1}=\frac{m}{2}$ and $s_{2}=\frac{d}{2}$. Right: $s_{1}=\frac{m}{3}$ and $s_{2}=\frac{d}{2}$.
		\end{tabular}
		\caption{Numerical experiments for MF  problem \eqref{SSNMF}. 
			%Numerical experiments for solving nonconvex sparse NMF problem \eqref{SSNMF}   on \emph{ORL} dataset with different $s_{1}$ and $s_{2}$. 
			%Top row: fixed $s_{1}$ and different $s_{2}$. Button row: different $s_{1}$ and fixed $s_{2}$.
		}
		\label{SSNMF_epochs}
	\end{figure}
	
	% The \emph{ORL} dataset with a reduced dimension of $r=25$ is utilized to demonstrate the efficiency of the proposed algorithm on the optimization problem \eqref{SSNMF}. All algorithms are run for $200$ epochs, with a subsampling rate of $5\%$ for stochastic algorithms. Figure \ref{SSNMF_epochs} compares the mean of the log of the objective function $\Phi(U_{k}, V_{k})$ versus the number of epochs for various algorithms, including BPG, BPGE, BPSG-SGD/SAGA/SARAH, and BPSGE-SGD/SAGA/SARAH.  
	%The first row of the figure displays the results with different values of $s_{1}$ and a fixed $s_{2}$, while the second row shows the results with a fixed $s_{1}$ and different values of $s_{2}$.

	% Fig.~\ref{SSNMF_epochs} shows that the BPGE algorithm outperforms the BPG algorithm, thanks to its extrapolation technique. The stochastic gradient algorithms that use the SAGA/SARAH estimator (BPSGE-SAGA/SARAH) perform better than their deterministic counterparts and the SGD estimator (BPSGE-SGD). The variance reduction technique also proves to be effective in the stochastic framework. Overall, the extrapolation technique speeds up convergence and the BPSGE-SAGA/SARAH algorithm provides the best results.

	Fig.~\ref{SSNMF_epochs} shows that the BPGE algorithm outperforms the BPG algorithm, thanks to its extrapolation technique. The adaptive stepsize and the variance reduction techniques also prove to be effective in the stochastic framework. 
	%The stochastic gradient algorithms that use the SAGA/SARAH estimator (BPSGE-SAGA/SARAH) perform better than their deterministic counterparts and the SGD estimator (BPSGE-SGD). The variance reduction technique also proves to be effective in the stochastic framework. Overall, the extrapolation technique speeds up convergence and the BPSGE-SAGA/SARAH algorithm provides the best results.

	%%%%%%%%%%%%%%%%%%%%%%%

	%%%%%%%%%%%%%%%%%%%
	\section{Conclusion}\label{conclusion}
	This paper presented a Bregman proximal stochastic gradient descent algorithm with extrapolation (BPSGE) for the objective function that  
	%characterized by a non-differentiable and a differentiable nonconvex function, where the differentiable component 
	lacks a global Lipschitz continuous gradient.  Under certain suitable conditions, we established the subsequential convergence, demonstrated that the subgradient of the objective function exhibits a sublinear convergence rate, and established the global convergence of the sequence. At last, we conducted numerical experiments on three specific applications and demonstrated the superior performance of the BPSGE algorithm. 
	
	% This paper presented a Bregman proximal stochastic gradient descent algorithm with extrapolation (BPSGE) for objective functions characterized by a non-differentiable and a differentiable nonconvex function, where the differentiable component lacks a global Lipschitz continuous gradient.  Under certain suitable conditions, we established the subsequential convergence of this algorithm and demonstrated that the subgradient of the objective function exhibits a sublinear convergence rate. Additionally, we also established the global convergence of the sequence. To evaluate the performance of the proposed algorithm, we conducted numerical experiments on three specific applications. These experiments demonstrated the superior performance of the BPSGE algorithm. 
	
	%in comparison to the existing approaches.
	\section{Acknowledgments}
	The authors are grateful to the anonymous referees for their valuable comments and suggestions.  This research is supported by the National Natural Science Foundation of China (NSFC) grants 12126608,  12131004, 12126603,  the R\&D project of Pazhou Lab (Huangpu) (Grant no. 2023K0603), and the Fundamental Research Funds for the Central Universities (Grant No. YWF-22-T-204).
	%%%%%%%%%%%
	\bibliography{aaai24}

	\onecolumn
	\appendix
	\section{Mathematical Proofs} \label{proof_details}
	
	\subsection{Proof of Lemma \ref{lemma_Phi_kk1}}
	\begin{proof}
		From the convexity of $h(\cdot)+\frac{\alpha}{2}\|\cdot\|^{2}$  in Assumption \ref{assume_01}, we obtain  
		%from the sub-gradient inequality that 
		\[
		h(x_{k+1})+\frac{\alpha}{2}\|x_{k+1}\|^{2}+\langle \xi_{k+1}+\alpha x_{k+1},x_{k}-x_{k+1}\rangle \le h(x_{k})+\frac{\alpha}{2}\|x_{k}\|^{2},
		\]
		where $\xi_{k+1}\in \partial h(x_{k+1})$. By rearranging the inequality we obtain
		\begin{eqnarray}
			h(x_{k+1})-\frac{\alpha}{2}\|x_{k+1}-x_{k}\|^{2}+\langle \xi_{k+1},x_{k}-x_{k+1}\rangle\le h(x_{k}). \label{h_ineq_01}
		\end{eqnarray}
		From the first-order optimality condition of  \eqref{xk_update} in Algorithm \ref{BPSGE}, it shows that
		\[
		\xi_{k+1} +\tilde{\nabla} f(\bar{x}_{k}) +\frac{1}{\eta_{k}}(\nabla \psi(x_{k+1})-\nabla\psi(\bar{x}_{k}))=0. 
		\]
		Combining the above inequality with \eqref{h_ineq_01}, we can derive the following result
		\begin{eqnarray}
			\begin{aligned}
				& h(x_{k+1})-\frac{\alpha}{2}\|x_{k+1}-x_{k}\|^{2}-\langle \tilde{\nabla} f(\bar{x}_{k}), x_{k}-x_{k+1}\rangle
				+\frac{1}{\eta_{k}}\langle \nabla \psi(\bar{x}_{k})-\nabla \psi(x_{k+1}), x_{k}-x_{k+1}\rangle\\
				=&h(x_{k+1})-\frac{\alpha}{2}\|x_{k+1}-x_{k}\|^{2}-\langle \tilde{\nabla} f(\bar{x}_{k}), x_{k}-x_{k+1}\rangle
				+\frac{1}{\eta_{k}}(D_{\psi}(x_{k},x_{k+1})+D_{\psi}(x_{k+1},\bar{x}_{k})-D_{\psi}(x_{k},\bar{x}_{k}))\\
				\le & h(x_{k}),
				\label{h_ineq_02}
			\end{aligned}
		\end{eqnarray}
		where the last equality follows from the three-point identity. 
		
		Furthermore, since $f$ is an $(\bar{L},\underline{L})$-relative smooth function with respect to $\psi$, we have
		\[
		f(x_{k+1}) \le f(\bar{x}_{k})+\langle \nabla f(\bar{x}_{k}), x_{k+1}-\bar{x}_{k}\rangle +\bar{L} D_{\psi}(x_{k+1},\bar{x}_{k}),
		\]
		and
		\[
		f(\bar{x}_{k})+\langle \nabla f(\bar{x}_{k}), x_{k}-\bar{x}_{k}\rangle \le f(x_{k}) +\underline{L} D_{\psi}(x_{k},\bar{x}_{k}).
		\]
		It shows that
		\begin{eqnarray}
			\quad\quad f(x_{k+1})\le f(x_{k})
			+\langle \nabla f(\bar{x}_{k}) ,x_{k+1}-x_{k}\rangle +\underline{L}D_{\psi}(x_{k},\bar{x}_{k})+\bar{L}D_{\psi}(x_{k+1},\bar{x}_{k}).\label{f_ineq_01}
		\end{eqnarray}

		By summing inequalities \eqref{h_ineq_02} and \eqref{f_ineq_01} together, we obtain
		\begin{eqnarray}
			\begin{aligned}
				\Phi(x_{k+1})\le &\Phi(x_{k}) +\frac{\alpha}{2}\|x_{k+1}-x_{k}\|^{2}+\langle \nabla f(\bar{x}_{k})-\tilde{\nabla} f(\bar{x}_{k}), x_{k+1}-x_{k}\rangle\\
				&+\left(\frac{1}{\eta_{k}}+\underline{L}\right)D_{\psi}(x_{k},\bar{x}_{k}) -\frac{1}{\eta_{k}}D_{\psi}(x_{k},x_{k+1})+\left(\bar{L}-\frac{1}{\eta_{k}}\right)D_{\psi}(x_{k+1},\bar{x}_{k})\\
				\le&\Phi(x_{k}) +\frac{\alpha+\gamma_{k}}{2}\|x_{k+1}-x_{k}\|^{2}+\frac{1}{2\gamma_{k}}\|\nabla f(\bar{x}_{k})-\tilde{\nabla} f(\bar{x}_{k})\|_{*}^{2}\\
				&+\left(\frac{1}{\eta_{k}}+\underline{L}\right)D_{\psi}(x_{k},\bar{x}_{k}) -\frac{1}{\eta_{k}}D_{\psi}(x_{k},x_{k+1}),
			\end{aligned}
		\end{eqnarray}
		where the last inequality follows from $\langle a,b\rangle \le \frac{\gamma}{2}\|a\|^{2}+\frac{1}{2\gamma}\|b\|^{2}$ for any $\gamma_k>0$ and $\eta_{k}\le \bar{L}^{-1}$. 
		
		By applying the conditional expectation operator $\mathbb{E}_{k}$ to the above inequality and bounding the MSE term by \eqref{MSE_l22} in Definition \ref{vr_definition},   we have 
		\begin{eqnarray}
			\begin{aligned}
				\mathbb{E}_{k}[\Phi(x_{k+1})]\le& \Phi(x_{k}) +\frac{\alpha+\gamma_{k}}{2}\mathbb{E}_{k}[\|x_{k+1}-x_{k}\|^{2}]+\frac{1}{2\gamma_{k}}\mathbb{E}_{k}[\|\nabla f(\bar{x}_{k})-\tilde{\nabla} f(\bar{x}_{k})\|_{*}^{2}]\\
				&+\left(\frac{1}{\eta_{k}}+\underline{L}\right)D_{\psi}(x_{k},\bar{x}_{k}) -\frac{1}{\eta_{k}}\mathbb{E}_{k}[D_{\psi}(x_{k},x_{k+1})]\\
				\le&\Phi(x_{k}) +\frac{\alpha+\gamma_{k}}{2}\mathbb{E}_{k}[\|x_{k+1}-x_{k}\|^{2}]+\frac{1}{2\gamma_{k}}\Gamma_{k}+ \frac{V_{1}}{2\gamma_{k}}\|x_{k}-x_{k-1}\|^{2}\\
				&+\frac{V_{1}}{2\gamma_{k}} \|x_{k-1}-x_{k-2}\|^{2}+\left(\frac{1}{\eta_{k}}+\underline{L}\right)D_{\psi}(x_{k},\bar{x}_{k}) -\frac{1}{\eta_{k}}\mathbb{E}_{k}[D_{\psi}(x_{k},x_{k+1})]\\ 
				\le&\Phi(x_{k}) +\frac{\alpha+\gamma_{k}}{2}\mathbb{E}_{k}[\|x_{k+1}-x_{k}\|^{2}]+\frac{1}{2\gamma_{k}\tau}(\Gamma_{k}-\mathbb{E}_{k}[\Gamma_{k+1}])\\
				&+\left(\frac{V_{\Gamma}}{2\gamma_{k}\tau}+\frac{V_{1}}{2\gamma_{k}}\right)\|x_{k}-x_{k-1}\|^{2}+\left(\frac{V_{\Gamma}}{2\gamma_{k}\tau}+\frac{V_{1}}{2\gamma_{k}}\right)\|x_{k-1}-x_{k-2}\|^{2}\\
				&+\left(\frac{1}{\eta_{k}}+\underline{L}\right)D_{\psi}(x_{k},\bar{x}_{k}) -\frac{1}{\eta_{k}}\mathbb{E}_{k}[D_{\psi}(x_{k},x_{k+1})],
			\end{aligned}\label{ineq_002}
		\end{eqnarray}
		where the last inequality follows from \eqref{Gamma_k1_k} in Definition \ref{vr_definition}.  From \eqref{extra_ineq} and $\eta_{k}\le\min\{\eta_{k-1},\bar{L}^{-1}\}$, it shows that
		\begin{eqnarray}
			\begin{aligned}
				\left(\frac{1}{\eta_{k}}+\underline{L}\right)D_{\psi}(x_{k},\bar{x}_{k})
				\le \frac{\underline{L}\eta_{k}+1}{\eta_{k}}\frac{\delta-\epsilon}{1+\underline{L}\eta_{k-1}}D_{\psi}(x_{k-1},x_{k})\le\frac{\delta-\epsilon}{\eta_{k}}D_{\psi}(x_{k-1},x_{k}).  \label{ineq_003}
			\end{aligned}
		\end{eqnarray}
		Combining \eqref{ineq_002} with \eqref{ineq_003}, we can get 
		\[
		\begin{aligned}
			\mathbb{E}_{k}[\Phi(x_{k+1})]\le& \Phi(x_{k}) -\left(\frac{1}{\eta_{k}}-\alpha-\gamma_{k}\right)\mathbb{E}_{k}[D_{\psi}(x_{k},x_{k+1})]+\frac{1}{2\gamma_{k}\tau}(\Gamma_{k}-\mathbb{E}_{k}[\Gamma_{k+1}])\\
			+&\left(\frac{\delta-\epsilon}{\eta_{k}}+\frac{V_{\Gamma}}{\gamma_{k}\tau}+\frac{V_{1}}{\gamma_{k}}\right)D_{\psi}(x_{k-1},x_{k})+\left(\frac{V_{\Gamma}}{\gamma_{k}\tau}+\frac{V_{1}}{\gamma_{k}}\right)D_{\psi}(x_{k-2},x_{k-1}).
		\end{aligned}
		\]
		Therefore, the results can be obtained by  rearranging the above terms  with $\gamma_{k}=\sqrt{2(V_{\Gamma}/\tau+V_{1})}$. 
		This completes the proof.
	\end{proof}
	
	%%%%%%%%%%%%%%%%%%%%%%%%%%%%%
	\subsection{Proof of  Lemma \ref{lyapunov_descent}}
	\begin{proof}
		From Lemma \ref{lemma_Phi_kk1}, it shows that
		\begin{eqnarray}
			\begin{aligned}
				\eta_{k}(\Phi(x_{k})-\mathcal{V}(\Phi))
				\ge&\eta_{k}(\mathbb{E}_{k}[\Phi(x_{k+1}])-\mathcal{V}(\Phi)) +(1-\eta_{k}\alpha-\eta_{k}\gamma)\mathbb{E}_{k}[D_{\psi}(x_{k},x_{k+1})]\\
				&+\frac{\eta_{k}}{2\tau\gamma}(\mathbb{E}_{k}[\Gamma_{k+1}]-\Gamma_{k})
				-\left(\delta-\epsilon+\frac{\gamma\eta_{k}}{2}\right)D_{\psi}(x_{k-1},x_{k})
				-\frac{\gamma\eta_{k}}{2}D_{\psi}(x_{k-2},x_{k-1}). \label{inequality_01}
			\end{aligned}
		\end{eqnarray}
		Combining \eqref{inequality_01} with $\eta_{k}\le\eta_{k-1}$, we have
		\begin{align*}
			&\Psi_{k}-\mathbb{E}_{k}[\Psi_{k+1}]\\
			=&\eta_{k-1}(\Phi(x_{k}) -\mathcal{V}(\Phi))  +\left(1-\eta_{k-1}\alpha-\eta_{k-1}\gamma-\frac{\epsilon}{3}\right)D_{\psi}(x_{k-1},x_{k})-\frac{\eta_{k}}{2\tau\gamma}\mathbb{E}_{k}[\Gamma_{k+1}]\\
			&+\frac{\eta_{k-1}}{2\tau\gamma}\Gamma_{k}+\eta_{k-1}\left(\frac{\gamma}{2}+\frac{\epsilon}{3\eta_{k-1}}\right)D_{\psi}(x_{k-2},x_{k-1})
			-\eta_{k}(\mathbb{E}_{k}[\Phi(x_{k+1})] -\mathcal{V}(\Phi))\\
			&-\eta_{k}\left(\frac{\gamma}{2}+\frac{\epsilon}{3\eta_{k}}\right)D_{\psi}(x_{k-1},x_{k}) -(1-\eta_{k}\alpha-\eta_{k}\gamma-\frac{\epsilon}{3})\mathbb{E}_{k}[D_{\psi}(x_{k},x_{k+1})]\\
			\ge&\eta_{k}(\Phi(x_{k}) -\mathcal{V}(\Phi))  +\left(1-\eta_{k-1}\alpha-\eta_{k-1}\gamma-\frac{\epsilon}{3}\right)D_{\psi}(x_{k-1},x_{k})-\frac{\eta_{k}}{2\tau\gamma}\mathbb{E}_{k}[\Gamma_{k+1}]\\
			&+\frac{\eta_{k}}{2\tau\gamma}\Gamma_{k}+\eta_{k-1}\left(\frac{\gamma}{2}+\frac{\epsilon}{3\eta_{k-1}}\right)D_{\psi}(x_{k-2},x_{k-1})
			-\eta_{k}(\mathbb{E}_{k}[\Phi(x_{k+1})] -\mathcal{V}(\Phi))\\
			&-\eta_{k}\left(\frac{\gamma}{2}+\frac{\epsilon}{3\eta_{k}}\right)D_{\psi}(x_{k-1},x_{k}) -\left(1-\eta_{k}\alpha-\eta_{k}\gamma-\frac{\epsilon}{3}\right)\mathbb{E}_{k}[D_{\psi}(x_{k},x_{k+1})]\\
			\ge&\left(1-\delta-\eta_{k-1}\alpha-(\eta_{k-1}+\eta_{k})\gamma\right)D_{\psi}(x_{k-1},x_{k})
			+\frac{\epsilon}{3}(\mathbb{E}_{k}D_{\psi}(x_{k},x_{k+1})+D_{\psi}(x_{k-1},x_{k})+D_{\psi}(x_{k-2},x_{k-1}))\\
			\ge&\left(1-\delta-\eta_{k-1}\alpha-2\eta_{k-1}\gamma\right)D_{\psi}(x_{k-1},x_{k})
			+\frac{\epsilon}{3}(\mathbb{E}_{k}D_{\psi}(x_{k},x_{k+1})+D_{\psi}(x_{k-1},x_{k})+D_{\psi}(x_{k-2},x_{k-1}))\\
			\ge&\frac{\epsilon}{3}(\mathbb{E}_{k}D_{\psi}(x_{k},x_{k+1})+D_{\psi}(x_{k-1},x_{k})+D_{\psi}(x_{k-2},x_{k-1})),
		\end{align*}
		where the second   and the last inequality follow from \eqref{inequality_01} and \eqref{stepsize_set}, respectively. This completes the proof. 
	\end{proof}
	%%%%%%%%%%%%%%%%%%%%%%%%%%%%
	\subsection{Proof of  Theorem \ref{subsequence_convergence}}
	\begin{proof}
		\begin{itemize}
			\item[(i)] This statement  follows   directly from Lemma \ref{lyapunov_descent} and $\epsilon>0$.
			\item[(ii)] By  summing \eqref{decent_inequality_01} from $k=0$ to a positive integer $K$, we have
			\[
			\sum_{k=1}^{K}\mathbb{E}[D_{\psi}(x_{k-1},x_{k})]\le \frac{3}{\epsilon}\mathbb{E}[\Psi_{1}-\Psi_{K+1}]\le \frac{3}{\epsilon}\Psi_{1},
			\]
			where the last inequality follows from $\Psi_{k}\ge 0$ for any $k>0$. Taking the limit as $K\rightarrow+\infty$, we have $\sum_{k=1}^{+\infty}\mathbb{E}[D_{\psi}(x_{k-1},x_{k})]<+\infty$. Then we may deduce that  the sequence $\{\mathbb{E}[ D_{\psi}(x_{k-1}, x_{k})]\}$ converges to zero. 
			\item[(iii)] We have
			\[
			K\min_{1\le k\le K}\mathbb{E}[D_{\psi}(x_{k-1},x_{k})]\le\sum_{k=1}^{K}\mathbb{E}[D_{\psi}(x_{k-1},x_{k})]\le\frac{3}{\epsilon}\Psi_{1},
			\]
			which yields the desired result. 
		\end{itemize}
		This completes the proof.
	\end{proof}

	%%%%%%%%%%%%%%%%%%%%%%%%%%%%%%%%%%%%%%%%%%%%%%%%%%%%%%%%%%%%%%
	\subsection{Proof of  Proposition \ref{vr_gra}}\label{proof_vr_gra}
	\begin{proof}
		For the SARAH stochastic gradient estimator,  we can get the results directly similar to the proof of Lemma 5 in \cite{WangH23}. 
		%The proof of   Proposition \ref{vr_gra}\,(i) is completed.  
		
		Now we  consider the proof of Proposition \ref{vr_gra}\,(ii). Firstly, we define the SAGA stochastic gradient estimator $\tilde{\nabla}^{SAGA}f(\bar{x}_{k})$ as 
		\[
		\tilde{\nabla}^{SAGA}f(\bar{x}_{k}):= \frac{1}{b}\sum_{j\in B_{k}}\left(\nabla f_{j}(\bar{x}_{k})-\nabla f_{j}(z_{k}^{j})\right)+\frac{1}{n}\sum_{i=1}^{n}\nabla f_{i}(z_{k}^{i}),
		\]
		where $ z_{k+1}^{i}=\begin{cases}
			\bar{x}_{k},\quad&\text{if }i\in B_{k},\\
			z_{k}^{i},\quad &\text{otherwise}.
		\end{cases}$ 
		From the Lipschitz continuity of $\nabla f_{i}(\cdot)$, it shows that
		\begin{align*}
			\mathbb{E}_{k}\|\tilde{\nabla}^{SAGA}f(\bar{x}_{k})-\nabla f(\bar{x}_{k})\|_{*}^{2}
			=&\mathbb{E}_{k}\|\frac{1}{b}\sum_{j\in B_{k}}\left(\nabla f_{j}(\bar{x}_{k})-\nabla f_{j}(z_{k}^{j})\right)+\frac{1}{n}\sum_{i=1}^{n}\nabla f_{i}(z_{k}^{i})-\nabla f(\bar{x}_{k})\|_{*}^{2}\\
			\le&\frac{1}{b^{2}}\sum_{j\in B_{k}}\|\nabla f_{j}(\bar{x}_{k})-\nabla f_{j}(z_{k}^{j})\|_{*}^{2}\\
			=&\frac{1}{bn}\sum_{i=1}^{n}\|\nabla f_{i}(\bar{x}_{k})-\nabla f_{i}(z_{k}^{i})\|_{*}^{2},
		\end{align*}
		where last inequality follows from the fact that $\mathbb{E}_{k}\|y_{1}+\cdots+y_{t}\|_{*}^{2}=\mathbb{E}_{k}\|y_{1}\|_{*}^{2}+\cdots+\mathbb{E}_{k}\|y_{t}\|_{*}^{2}$ for any independent random variables $y_{i} (i=1,\dots,t)$ with $\mathbb{E}_{k}[y_{i}]=0$ for all $i$. Combined with Jensen's inequality, we can get 
		\begin{align*}
			\mathbb{E}_{k}\|\tilde{\nabla}^{SAGA}f(\bar{x}_{k})-\nabla f(\bar{x}_{k})\|_{*}
			\le&\sqrt{\mathbb{E}_{k}\|\tilde{\nabla}^{SAGA}f(\bar{x}_{k})-\nabla f(\bar{x}_{k})\|_{*}^{2}}\\
			\le&\frac{1}{\sqrt{bn}}\sqrt{\sum_{i=1}^{n}\|\nabla f_{i}(\bar{x}_{k})-\nabla f_{i}(z_{k}^{i})\|_{*}^{2}}\\
			\le&\frac{1}{\sqrt{bn}}\sum_{i=1}^{n}\|\nabla f_{i}(\bar{x}_{k})-\nabla f_{i}(z_{k}^{i})\|_{*}.
		\end{align*}
		We bound the MSE of the stochastic gradient estimator $\tilde{\nabla}^{SAGA}f(\cdot)$ as follows,
		\begin{align*}
			&\frac{1}{bn}\sum_{i=1}^{n}\mathbb{E}_{k}\|\nabla f_{i}(\bar{x}_{k})-\nabla f_{i}(z_{k}^{i})\|_{*}^{2}\\
			\le&\frac{1+\delta}{bn}\mathbb{E}_{k}\sum_{i=1}^{n}\|\nabla f_{i}(\bar{x}_{k-1})-\nabla f_{i}(z_{k}^{i})\|_{*}^{2}+\frac{1+\delta^{-1}}{bn}\sum_{i=1}^{n}\|\nabla f_{i}(\bar{x}_{k})-\nabla f_{i}(\bar{x}_{k-1})\|_{*}^{2}\\
			\le&\frac{1+\delta}{bn}(1-\frac{b}{n})\sum_{i=1}^{n}\|\nabla f_{i}(\bar{x}_{k-1})-\nabla f_{i}(z_{k-1}^{i})\|_{*}^{2}+\frac{1+\delta^{-1}}{b}M_{1}^{2}\|\bar{x}_{k}-\bar{x}_{k-1}\|^{2}\\
			\le&\frac{1+\delta}{bn}(1-\frac{b}{n})\sum_{i=1}^{n}\|\nabla f_{i}(\bar{x}_{k-1})-\nabla f_{i}(z_{k-1}^{i})\|_{*}^{2}+\frac{1+\delta^{-1}}{b}M_{1}^{2}[(1+\beta_{k}^{2})\|x_{k}-x_{k-1}\|^{2}
			+\beta_{k-1}^{2}\|x_{k-1}-x_{k-2}\|^{2}]\\
			\le&\frac{1+\delta}{bn}(1-\frac{b}{n})\sum_{i=1}^{n}\|\nabla f_{i}(\bar{x}_{k-1})-\nabla f_{i}(z_{k-1}^{i})\|_{*}^{2}+\frac{2+2\delta^{-1}}{b}M_{1}^{2}[\|x_{k}-x_{k-1}\|^{2}
			+\|x_{k-1}-x_{k-2}\|^{2}],
		\end{align*}
		where the first inequality follows from $\|x-z\|_{*}^{2}\le (1+\delta)\|x-y\|_{*}^{2}+(1+\delta^{-1})\|y-z\|_{*}^{2}$. 
		Let $\Gamma_{k+1}:=\frac{1}{bn}\sum_{i=1}^{n}\|\nabla f_{i}(\bar{x}_{k})-\nabla f_{i}(z_{k}^{i})\|_{*}^{2}$ and $\delta=\frac{b}{2n}$, it shows that 
		\[
		\begin{aligned}
			\mathbb{E}_{k}\Gamma_{k+1}\le& (1+\frac{b}{2n})(1-\frac{b}{n})\Gamma_{k}+\frac{2+\frac{4n}{b}}{b}M_{1}^{2}[\|x_{k}-x_{k-1}\|^{2}+\|x_{k-1}-x_{k-2}\|^{2}]\\
			\le&(1-\frac{b}{2n})\Gamma_{k}+\frac{2b+4n}{b^{2}} / 2n+\frac{4n^{2}}{b}M_{1}^{2}[\|x_{k}-x_{k-1}\|^{2}+\|x_{k-1}-x_{k-2}\|^{2}]. 
		\end{aligned}
		\] 
		This proves the geometric decay of $\Gamma_{k}$ in expectation. Similar to Appendix B in \cite{DriggsTLDS2020}, we also have that the third condition holds in Definition \ref{vr_definition}. This completes the proof.
	\end{proof}
	
	%%%%%%%%%%%%%%%%%%%
	\subsection{Proof of  Lemma \ref{subgradient_bound}}
	\begin{proof}
		From the implicit deﬁnition of the proximal operator \eqref{xk_update} in the BPSGE algorithm, we have that
		\[
		0\in \partial h(x_{k+1}) +\tilde{\nabla}f(\bar{x}_{k})+\frac{1}{\eta_{k}}(\nabla \psi(x_{k+1})-\nabla \psi(\bar{x}_{k})).
		\]
		Combining it with  $\partial \Phi(x_{k+1})\equiv\nabla f(x_{k+1})+\partial h(x_{k+1})$,   we have $w_{k+1}\in \partial \Phi(x_{k+1})$. All that remains is to bound the norm of $w_{k+1}$.  $\nabla f$ and $\nabla \psi$ are Lipschitz continuous  with constants $M_{1}$ and $M_{2}$ on any bounded subset of $\mathbb{R}^{d}$, respectively (See Assumption \ref{assume_03}). It shows that
		\begin{align*}
			&\mathbb{E}_{k}\|w_{k+1}\|\\
			\le&\mathbb{E}_{k}\|\nabla f(x_{k+1})-\tilde{\nabla}f(\bar{x}_{k})+\frac{1}{\eta_{k}}(\nabla \psi(\bar{x}_{k})-\nabla \psi(x_{k+1}))\|\\
			\le&\mathbb{E}_{k}\|\nabla f(x_{k+1})-\tilde{\nabla}f(\bar{x}_{k})\|+\frac{1}{\eta_{k}}\mathbb{E}_{k}\|\nabla \psi(\bar{x}_{k})-\nabla \psi(x_{k+1})\|\\
			\le&\mathbb{E}_{k}\|\nabla f(x_{k+1})-\nabla f(\bar{x}_{k})\|+\mathbb{E}_{k}\|\nabla f(\bar{x}_{k})-\tilde{\nabla} f(\bar{x}_{k})\|+\frac{1}{\eta_{k}}\mathbb{E}_{k}\|\nabla \psi(\bar{x}_{k})-\nabla \psi(x_{k+1})\|\\
			\le&M_{1}\mathbb{E}_{k}\|x_{k+1}-\bar{x}_{k}\|+\Upsilon_{k}+V_{2}\|x_{k}-x_{k-1}\|+V_{2}\|x_{k-1}-x_{k-2}\|+\frac{M_{2}}{\eta_{k}}\mathbb{E}_{k}\|x_{k+1}-\bar{x}_{k}\|\\
			\le& \left(M_{1}+\frac{M_{2}}{\eta_{k}}\right)\mathbb{E}_{k}\|x_{k+1}-x_{k}\|+\left(V_{2}+\beta_{k}M_{1}+\frac{\beta_{k}M_{2}}{\eta_{k}}\right)\|x_{k}-x_{k-1}\|
			+V_{2}\|x_{k-1}-x_{k-2}\|+\Upsilon_{k}\\
			\le&\left(M_{1}+\frac{M_{2}}{\eta}\right)\mathbb{E}_{k}\|x_{k+1}-x_{k}\|+\left(V_{2}+\beta_{k}M_{1}+\frac{\beta_{k}M_{2}}{\eta}\right)\|x_{k}-x_{k-1}\|
			+V_{2}\|x_{k-1}-x_{k-2}\|+\Upsilon_{k}\\
			\le& \rho(\mathbb{E}_{k}\|x_{k+1}-x_{k}\|+\|x_{k}-x_{k-1}\|+\|x_{k-1}-x_{k-2}\|) +\Upsilon_{k},
		\end{align*}
		where $\rho=\max\left\{M_{1}+\frac{M_{2}}{\eta}, V_{2}+\beta_{k}M_{1}+\frac{\beta_{k}M_{2}}{\eta}, V_{2}\right\}$. This completes the proof.
	\end{proof}
	%%%%%%%%%%%%%%%%%%%%%%%%%%%%%%%%%%%%%%%%%
	\subsection{Proof of  Lemma \ref{lemma_dist2}}
	\begin{proof}
		From Lemma \ref{subgradient_bound}, it shows that
		\begin{align*}
			&\mathbb{E}_{k}\|w_{k+1}\|^{2}\\
			\le&3\mathbb{E}_{k}\|\nabla f(x_{k+1})-\nabla f(\bar{x}_{k})\|^{2}+3\mathbb{E}_{k}\|\nabla f(\bar{x}_{k})-\tilde{\nabla} f(\bar{x}_{k})\|^{2}
			+\frac{3}{\eta_{k}}\mathbb{E}_{k}\|\nabla \psi(\bar{x}_{k})-\nabla \psi(x_{k+1})\|^{2}\\
			\le&3M_{1}^{2}\mathbb{E}_{k}\|x_{k+1}-\bar{x}_{k}\|^2+3\Gamma_{k}+3V_{1}\|x_{k}-x_{k-1}\|^{2}+3V_{1}\|x_{k-1}-x_{k-2}\|^{2}
			+\frac{3M_{2}^{2}}{\eta_{k}}\mathbb{E}_{k}\|x_{k+1}-\bar{x}_{k}\|^{2}\\
			\le& \left(6M_{1}^{2}+\frac{6M_{2}^{2}}{\eta_{k}}\right)\mathbb{E}_{k}\|x_{k+1}-x_{k}\|^{2}+\left(3V_{1}+6\beta_{k}^{2}M_{1}^{2}+\frac{6\beta_{k}^{2}M_{2}^{2}}{\eta_{k}}\right)\|x_{k}-x_{k-1}\|^{2}
			+3V_{1}\|x_{k-1}-x_{k-2}\|^{2}+3\Gamma_{k}\\
			\le&\left(6M_{1}^{2}+\frac{6M_{2}^{2}}{\eta}\right)\mathbb{E}_{k}\|x_{k+1}-x_{k}\|^{2}+\left(3V_{1}+6\beta_{k}^{2}M_{1}^{2}+\frac{6\beta_{k}^{2}M_{2}^{2}}{\eta}\right)\|x_{k}-x_{k-1}\|^{2}
			+3V_{1}\|x_{k-1}-x_{k-2}\|^{2}+3\Gamma_{k}\\
			\le& \bar{\rho}(\mathbb{E}_{k}\|x_{k+1}-x_{k}\|^{2}+\|x_{k}-x_{k-1}\|^{2}+\|x_{k-1}-x_{k-2}\|^{2}) +3\Gamma_{k},
		\end{align*}
		where $\bar{\rho}:=\max\left\{6M_{1}^{2}+\frac{6M_{2}^{2}}{\eta}, 3V_{1}+6\beta_{k}^{2}M_{1}^{2}+\frac{6\beta_{k}^{2}M_{2}^{2}}{\eta}, 3V_{1}\right\}$. By $\mathrm{dist}(0,\partial\Phi(x_{k+1}))^{2}\le\|w_{k+1}\|^{2}$ and taking full expectation on both sides, it shows that
		\[
		\mathbb{E}[\mathrm{dist}(0,\partial\Phi(x_{k+1}))^{2}]\le \bar{\rho}\mathbb{E}[\|x_{k+1}-x_{k}\|^{2}+\|x_{k}-x_{k-1}\|^{2}+\|x_{k-1}-x_{k-2}\|^{2}] +3\mathbb{E}\Gamma_{k}.
		\]
		This completes the proof.
	\end{proof}
	
	%%%%%%%%%%%%%%%%%%%%%%%%%%%%%
	\subsection{Proof of Theorem \ref{subgradient_rate}}
	\begin{proof}
		From Corollary \ref{lemma_Phi_kk1_Rd} and Corollary \ref{lemma_dist2}, it shows that
		\begin{align*}
			&\mathbb{E}[\Psi_{k}-\Psi_{k+1}]\\
			\ge&\frac{\epsilon}{6}\mathbb{E}[\|x_{k+1}-x_{k}\|^{2}+\|x_{k}-x_{k-1}\|^{2}+\|x_{k-1}-x_{k-2}\|^{2}]\\
			\ge&\sigma\mathbb{E}[\|x_{k+1}-x_{k}\|^{2}+\|x_{k}-x_{k-1}\|^{2}+\|x_{k-1}-x_{k-2}\|^{2}]
			+\frac{\epsilon/6-\sigma}{\bar{\rho}}\mathbb{E}[\mathrm{dist}(0,\partial\Phi(x_{k+1}))^{2}]-\frac{\epsilon/2-3\sigma}{\bar{\rho}}\mathbb{E}\Gamma_{k}\\
			\ge&\sigma\mathbb{E}[\|x_{k+1}-x_{k}\|^{2}+\|x_{k}-x_{k-1}\|^{2}+\|x_{k-1}-x_{k-2}\|^{2}]+\frac{\epsilon/2-3\sigma}{\tau\bar{\rho}}\mathbb{E}[\Gamma_{k+1}-\Gamma_{k}]\\
			&-\frac{(\epsilon/2-3\sigma)V_{\Gamma}}{\tau\bar{\rho}}\mathbb{E}[\|x_{k}-x_{k-1}\|^{2}+\|x_{k-1}-x_{k-2}\|^{2}]
			+\frac{\epsilon/6-\sigma}{\bar{\rho}}\mathbb{E}[\mathrm{dist}(0,\partial\Phi(x_{k+1}))^{2}]\\
			\ge&\sigma\mathbb{E}[\|x_{k+1}-x_{k}\|^{2}+\|x_{k}-x_{k-1}\|^{2}+\|x_{k-1}-x_{k-2}\|^{2}]
			+\frac{\epsilon/2-3\sigma}{\tau\bar{\rho}}\mathbb{E}[\Gamma_{k+1}-\Gamma_{k}]\\
			& -\frac{(\epsilon/2-3\sigma)V_{\Gamma}}{\tau\bar{\rho}}\mathbb{E}[\|x_{k+1}-x_{k}\|^{2}+\|x_{k}-x_{k-1}\|^{2}+\|x_{k-1}-x_{k-2}\|^{2}]
			+\frac{\epsilon/6-\sigma}{\bar{\rho}}\mathbb{E}[\mathrm{dist}(0,\partial\Phi(x_{k+1}))^{2}],
		\end{align*}
		where the third inequality follows from \eqref{Gamma_k1_k} in Definition \ref{vr_definition}. If we let $\sigma=\frac{(\epsilon/2-3\sigma)V_{\Gamma}}{\tau\bar{\rho}}$, i.e., $\sigma=\frac{\frac{\epsilon}{2}V_{\Gamma}}{3V_{\Gamma}+\tau\bar{\rho}}$, it shows that
		\[
		\mathbb{E}[\Psi_{k}-\Psi_{k+1}]\ge \frac{\epsilon/6-\sigma}{\bar{\rho}}\mathbb{E}[\mathrm{dist}(0,\partial\Phi(x_{k+1}))^{2}]+\frac{\epsilon/2-3\sigma}{\tau\bar{\rho}}\mathbb{E}[\Gamma_{k+1}-\Gamma_{k}].
		\]
		Summing up $k=1$ to $K$, we have
		\[
		\mathbb{E}[\Psi_{1}-\Psi_{K+1}]\ge \frac{\epsilon/6-\sigma}{\bar{\rho}}\sum_{k=1}^{K}\mathbb{E}[\mathrm{dist}(0,\partial\Phi(x_{k+1}))^{2}]+\frac{\epsilon/2-3\sigma}{\tau\bar{\rho}}\mathbb{E}[\Gamma_{K+1}-\Gamma_{1}],
		\]
		which means  there exists a  $k^{\prime} \in\{2,\dots,K+1\}$ such that
		\begin{align*}
			\mathbb{E}[\mathrm{dist}(0,\partial\Phi(x_{k^{\prime}}))^{2}]\le&\frac{1}{K}\sum_{k=1}^{K}\mathbb{E}[\mathrm{dist}(0,\partial\Phi(x_{k+1}))^{2}]\\
			\le&\frac{\bar{\rho}}{(\epsilon/6-\sigma)K}(\mathbb{E}[\Psi_{1}-\Psi_{K+1}]+\frac{\epsilon/2-3\sigma}{\tau\bar{\rho}}\mathbb{E}[\Gamma_{1}-\Gamma_{K+1}])\\
			\le&\frac{\bar{\rho}}{(\epsilon/6-\sigma)K}(\mathbb{E}\Psi_{1}+\frac{\epsilon/2-3\sigma}{\tau\bar{\rho}}\mathbb{E}\Gamma_{1}).
		\end{align*}
		This completes the proof.
	\end{proof}
	%%%%%%%%%%%%%%%%%%%%%%%%%%%%%%%%%%%
	
	\subsection{Proof of  Lemma \ref{statements_lemma}}
	\begin{proof}
		The four statements can be easily obtained by Corollary \ref{lemma_Phi_kk1_Rd} and Lemma \ref{subgradient_bound}, so we omit the details here for simplicity.
	\end{proof}
	The following lemma is from \cite{DriggsTLDS2020}, which is analogous to the Uniformized K{\L} property of \cite{BolteST14} and allows us to apply the K{\L} inequality.
	\begin{lemma}\label{E_KL_inequality}
		Let $\{x_{k}\}_{k=0}^{\infty}$ be a bounded sequence of iterates generated by the BPSGE algorithm using a variance-reduced gradient estimator (see Definition \ref{vr_definition}), and let $\Phi$ be a semialgebraic function satisfying the K{\L} property \cite{BolteST14} with exponent $\theta$. Then there exists an index $\bar{k}$ and a desingularizing function $\phi(r) = ar^{1-\theta}$ with $a>0$, $\theta\in[0,1)$  so that the following bound holds almost surely (a.s.),
		\begin{eqnarray}
			\phi'(\mathbb{E}[\Phi(x_{k})-\Phi_{k}^{*}])\mathbb{E}\mbox{dist}(0,\partial\Phi(x_{k}))\ge 1,\,\,\forall k>\bar{k},
		\end{eqnarray}
		where $\Phi_{k}^{*}$ is a nondecreasing sequence  converging to $\mathbb{E}\Phi(x_{*})$ for some $x_{*}\in\omega(x_{0})$.
	\end{lemma}
	
	\subsection{Proof of  Theorem \ref{global_convergence}}   
	\begin{proof}
		According to Lemma \ref{E_KL_inequality}, if $\Phi$ is a proper, lower semi-continuous, and semi-algebraic function, it will satisfy the
		KL property at any point of $\text{dom }\Phi$. Under Lemma \ref{E_KL_inequality}, combining Corollary \ref{lemma_Phi_kk1_Rd} with Lemma \ref{subgradient_bound}, we can get that the generated sequence
		$\{x_{k}\}$ is a Cauchy sequence which yields the result. The detailed proof of this theorem is similar to Theorem 2 in \cite{WangH23}. Thus the details are omitted here.
	\end{proof} 
	
	\subsection{Proof of Proposition \ref{Prop_GNMF}}
	\begin{proof}
		Combining Proposition C.4 with Proposition D.1 in \cite{MukkamalaO19}, we can directly get this result. 
	\end{proof}
	
	\subsection{Proof of  Proposition \ref{prop_WCMF}}
	\begin{proof}
		From the  update step for solving \eqref{xk_update}, we have  
		\begin{align*}
			&(U_{k+1},V_{k+1})\\
			\in&\,\,\underset{U\in\mathbb{R}^{m\times r},V\in\mathbb{R}^{r\times d}}{\arg\min}\,\,\left\{\eta_{k} h(U,V)+\langle P_{k},U\rangle+\langle Q_{k},V\rangle+\psi(U,V)\right\}\\
			=&\underset{U\in\mathbb{R}^{m\times r},V\in\mathbb{R}^{r\times d}}{\arg\min}\,\,\Big\{\eta_{k}\left(\lambda_{1}\|U\|_{1}-\frac{\lambda_{2}}{2}\|U\|_{F}^{2}\right)+\langle P_{k},U\rangle+\langle Q_{k},V\rangle
			+3\psi_{1}(U,V)+\|M\|_{F}\psi_{2}(U,V)+\frac{\eta_{k}\lambda_{2}}{2}\|U\|_{F}^{2}\Big\}\\
			=&\underset{U\in\mathbb{R}^{m\times r},V\in\mathbb{R}^{r\times d}}{\arg\min}\,\,\left\{\eta_{k}\|U\|_{1}+\langle P_{k},U\rangle+\langle Q_{k},V\rangle+3\psi_{1}(U,V)+\|M\|_{F}\psi_{2}(U,V)\right\}.
		\end{align*}
		From the Proposition C.5 in \cite{MukkamalaO19}, we can directly get the closed form of $(U_{k+1},V_{k+1})$.
	\end{proof}
	
	%%%%%%%%%%%%%%%%%%
	\section{More Details in the Numerical Experiments}\label{more_results}
	
	%\subsection{Numerical setup}
	
	We combine BPSGE with classic stochastic gradient estimator \cite{RobbinsM1951} (BPSGE-SGD), SAGA gradient estimator \cite{DefazioBL14} (BPSGE-SAGA), and SARAH gradient estimator \cite{NguyenLST17} (BPSGE-SARAH), and compare BPSGE-SGD/SAGA/SARAH with BPG \cite{BolteSTV18First},   BPGE (the special case of CoCaIn BPG \cite{MukkamalaOPS20}), and BPSG-SGD/SAGA/SARAH \cite{WangH23}. We consider three applications: graph regularized NMF, MF with weakly-convex regularization, and NMF with nonconvex sparsity constraints.

	Statistics of the four datasets in graph regularized NMF is listed in the following table. 
	\begin{table}[!ht]
		\centering
		%\setlength{\tabcolsep}{3pt}
		%\fontsize{9}{10}\selectfont
		\begin{tabular}{c|c|c|c}\hline
			Dataset& Size &Dimensionality     & Number of classes \\\hline
			\emph{COIL20} & 1440 & 1024 &20\\
			\emph{PIE} & 2856 &1024 & 68\\
			\emph{COIL100} & 7200 & 1024 &100\\
			\emph{TDT2} & 9394 & 36771 & 30\\\hline
		\end{tabular}
		\caption{Statistics of the four datasets in graph regularized NMF.}
		\label{GNMG_datasets}
	\end{table}

\begin{figure*}[!ht]
		\setlength\tabcolsep{1pt}
		\centering
		\begin{tabular}{cccc}
			\includegraphics[width=0.24\textwidth]{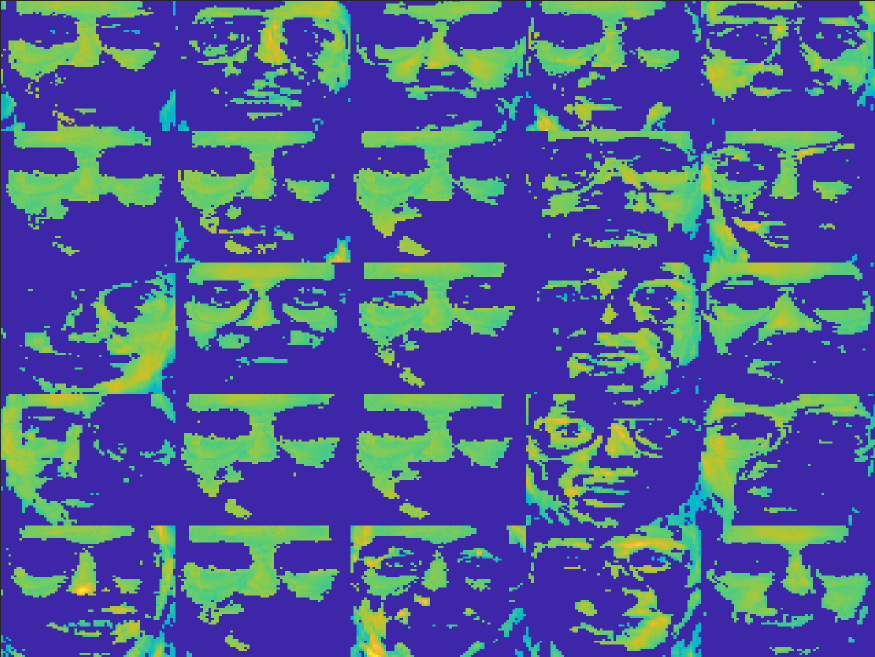}
			&\includegraphics[width=0.24\textwidth]{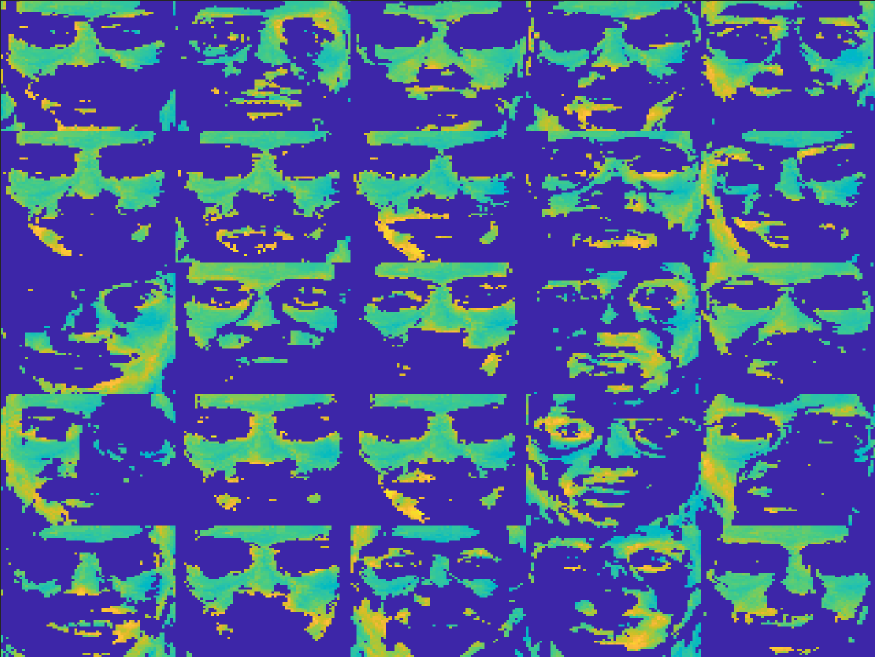}
			&\includegraphics[width=0.24\textwidth]{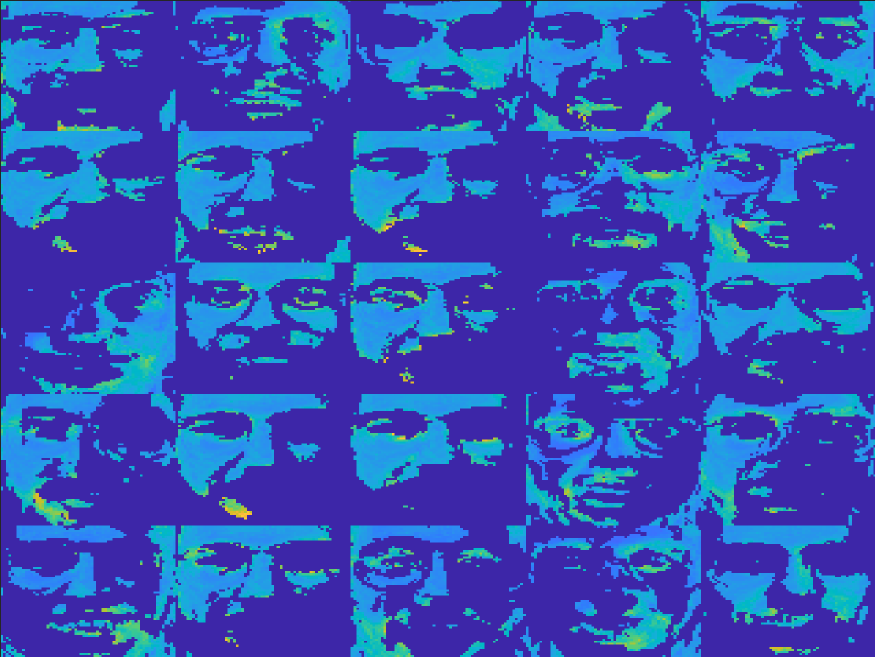}
			&\includegraphics[width=0.24\textwidth]{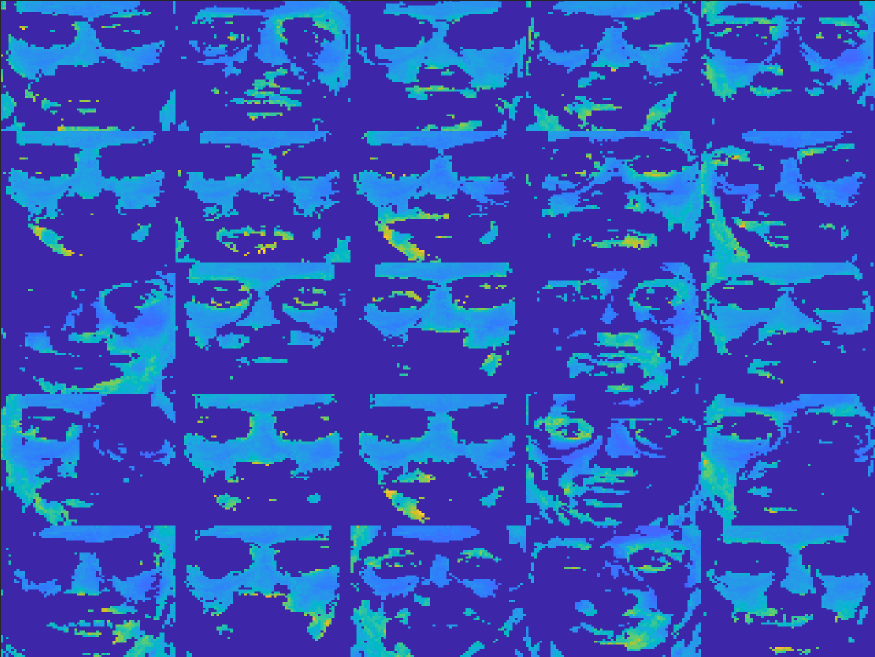}\\
			(a)BPG&(b)BPSG-SGD&(c)BPSG-SAGA&(d)BPSG-SARAH\\
			\includegraphics[width=0.24\textwidth]{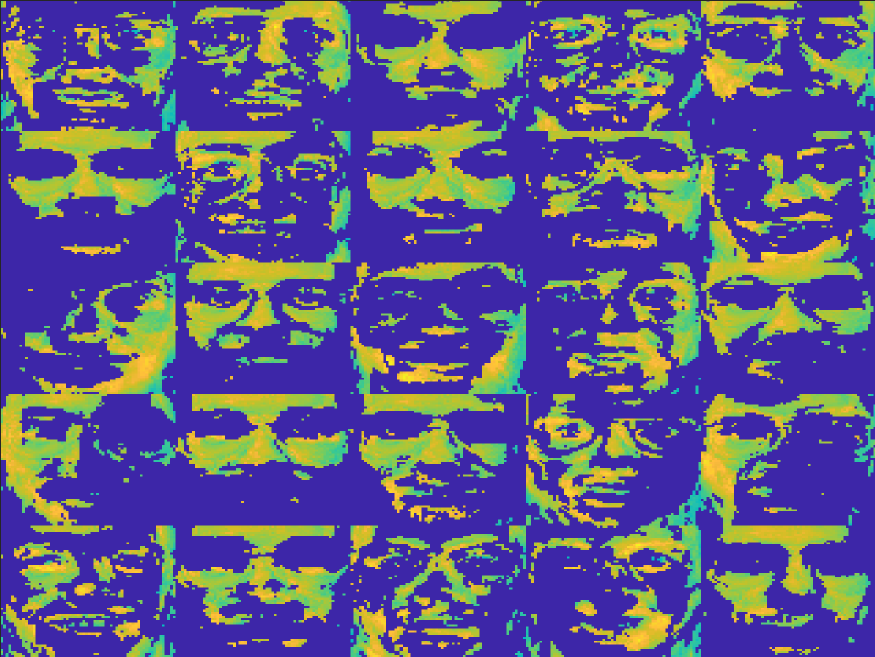}
			&\includegraphics[width=0.24\textwidth]{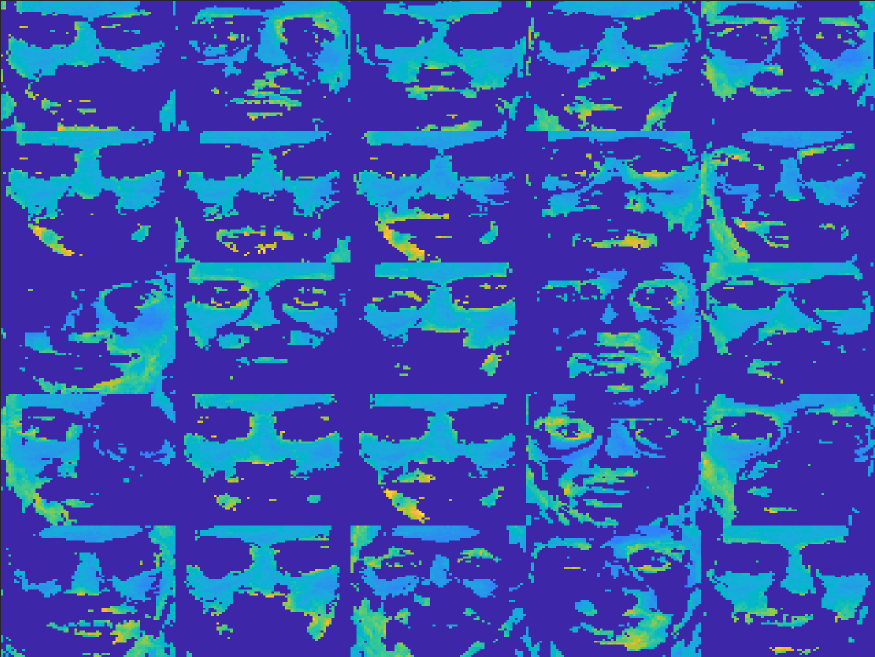}
			&\includegraphics[width=0.24\textwidth]{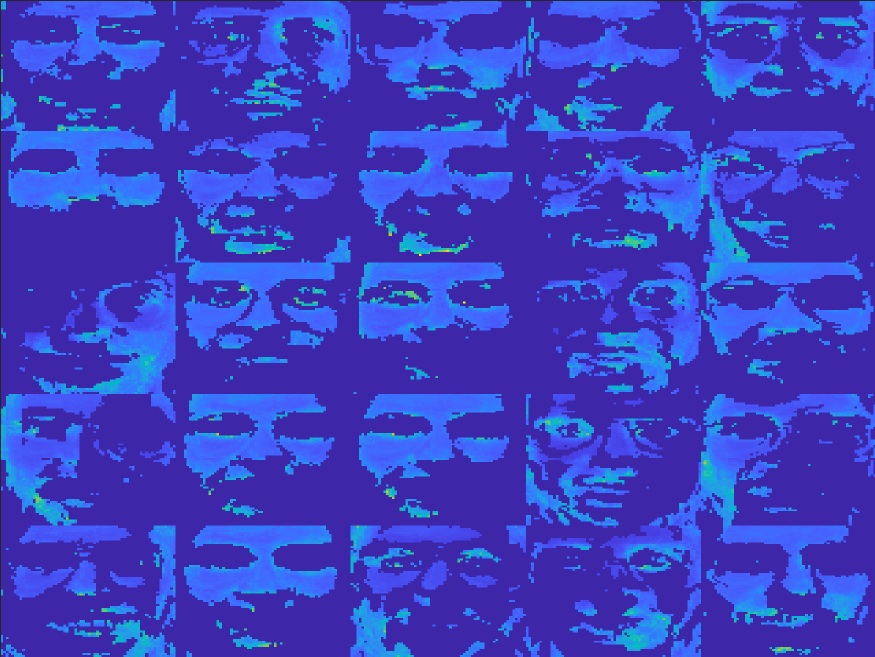}
			&\includegraphics[width=0.24\textwidth]{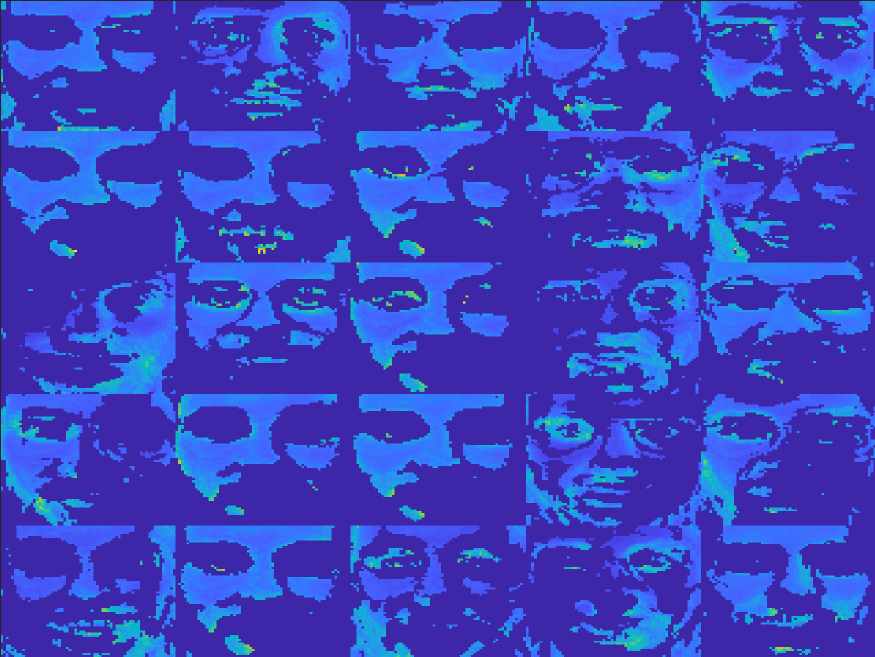}\\
			(e)BPGE&(f) BPSGE-SGD&(g)BPSGE-SAGA&(h)BPSGE-SARAH
		\end{tabular}
		\caption{The basis images generated by solving the nonconvex sparsity constrained NMF problem \eqref{SSNMF}  with $s_{1}=\frac{m}{3}$ and $s_{2}=\frac{d}{2}$.}
		\label{basis_SSNMF_ORL_32}
	\end{figure*}

	The basis images generated by solving the nonconvex sparsity constrained NMF problem \eqref{SSNMF} for $s_{1}=\frac{m}{3}$ and $s_{2}=\frac{d}{2}$ are shown in Figure \ref{basis_SSNMF_ORL_32}. It validates the BPSGE algorithm outperforms other determinant algorithms and stochastic algorithms without extrapolation. 
	
\end{document}